\documentclass{article}

\usepackage{arxiv}

\usepackage{natbib}
\usepackage{graphicx}
\usepackage{epstopdf, epsfig}
\usepackage{amssymb}
\usepackage{amsmath}
\usepackage{hyperref}
\usepackage[nameinlink, noabbrev]{cleveref}    
\usepackage{xcolor}
\usepackage{subfig}
\usepackage{float}
\usepackage[normalem]{ulem}
\usepackage{soul}
\usepackage{colortbl,xcolor}
\usepackage[ruled,linesnumbered]{algorithm2e}
\usepackage{multirow}
\usepackage{multicol}
\usepackage{tabularx}
\usepackage{setspace}

\hypersetup{
    colorlinks=true,
    linkcolor=blue,
    citecolor=blue,
    filecolor=magenta,      
    urlcolor=cyan,
    hypertexnames=true
    }

\usepackage{booktabs}       
\usepackage{nicefrac}       
\usepackage{microtype}      

\usepackage{ulem}

\title{An adaptive switch strategy for acquisition functions in Bayesian optimization of wind farm layout}

\author{Zhen-Fan Wang$^1$,
        Yu Tu$^1$,
        Kai Zhang$^1$\thanks{Corresponding author: kai.zhang@sjtu.edu.cn} ,
        Dai Zhou$^{1, 2}$,
        Onur Bilgen$^3$ \\
        $^1$School of Ocean and Civil Engineering, Shanghai Jiao Tong University, Shanghai, 200240, China \\
        $^2$Shenzhen Research Institute of Shanghai Jiao Tong University, Shenzhen, 518063, China \\
        $^3$Department of Mechanical and Aerospace Engineering, Rutgers University, New Brunswick, New Jersey 08901, USA
        }

\begin{document}

\maketitle
\begin{abstract}

Wind farm layout optimization, which seeks to minimizing wake losses and maximizing annual energy production by strategically adjusting wind turbines' location, is essential for the development of large-scale wind farms. While low-fidelity methods dominate wind farm layout optimization studies due to their computational efficiency, high-fidelity methods are less commonly applied due to their significant computational costs.
This paper introduces a Bayesian optimization framework that leverages a novel adaptive acquisition function switching strategy to enhance the efficiency and effectiveness of wind farm layout optimization using high-fidelity modeling methods. 
The proposed switch acquisition functions strategy alternates between maximum surrogate prediction and maximum-value entropy search acquisition functions, dynamically balancing exploration and exploitation during the optimization process. 
By iteratively retraining the Kriging model with intermediate optimal layouts, the framework progressively refines its predictions to accelerate convergence to optimal solutions.
The performance of the switch-acquisition-function-based Bayesian optimization framework is first validated using 4- and 10-dimensional Ackley benchmark functions, where it demonstrates superior optimization efficiency compared to conventional approaches using maximum surrogate prediction or maximum-value entropy search alone. 
The framework is then applied to wind farm layout optimization problems using computationally efficient Gaussian wake models for three wind farm cases with varying boundary constraints and wind distributions. 
Results show that the switch-acquisition-function-based Bayesian optimization framework outperforms traditional heuristic algorithms, achieving near-optimal annual energy output with significantly fewer calculations.
Finally, the framework is extended to high-fidelity wind farm layout optimization by coupling it with computational fluid dynamics  simulations, where turbine rotors are modeled as actuator disks. 
The novel switch-acquisition-function-based Bayesian optimization enables more effective exploration to achieve higher annual energy production in wind farm layout optimization, advancing the design of more effective wind farm layouts and contributing to the development of more efficient wind energy projects.

\end{abstract}

\section{Introduction}
\label{sec:intro}
Over the years, concerns regarding environmental degradation and global warming driven by the consumption of fossil fuels have intensified. The renewable energy sources, such as wind, solar, and tidal energy, have emerged as viable solutions to address these pressing issues \citep{spiru2024wind}. Among these, wind energy is widely regarded as one of the most promising renewable energy sources, owing to its abundant availability and significant potential.
As a result, an increasing number of countries have actively committed to the establishment of large-scale wind farms to exploit wind energy resources effectively. However, the efficiency of power generation in large wind farms is often compromised by wake interactions between turbines. These interactions result from the wake effects generated by upstream turbines, which cause reduced wind speed and increased turbulence intensity for downstream turbines.
The impact of such wake effects can be substantial, with average power losses reported to range between 10\% and 20\% in some large offshore wind farms \citep{barthelmie2009modelling}. Therefore, mitigating wake effects within wind farms is essential to enhance energy output efficiency and realize the full potential of wind energy.

As a promising means of alleviating the wake effects, wind farm layout optimization (WFLO) has been a subject of extensive research in recent decades \citep{erouglu2012design,stanley2019massive,antonini2020optimal,cruz2020wind,reddy2020wind,dong2021intelligent,thomas2023comparison}.
The goal of WFLO is to strategically arrange wind turbines within a given area to maximize the overall energy production of the wind farm.
A critical aspect of the optimization process is the adoption of reliable methods to evaluate the annual energy production (AEP) of different wind turbine layouts. 
Two primary approaches are commonly used in WFLO: low-fidelity methods based on analytical wake models, and high-fidelity methods employing computational fluid dynamics (CFD).
Analytical wake models such as the Jensen wake model \citep{jensen1983note} and the Gaussian model \citep{larsen1988simple} apply simplified analytical functions that depict the wake profiles based on reduced flow physics.
These low-fidelity methods are computationally efficient and are widely used in WFLO studies. 
However, the reliance on simplified physics often results in less accurate assessments of wake interactions, thereby leading to potentially erroneous AEP estimates.

Computational fluid dynamics (CFD) simulations offer more accurate insights into the wake development behind wind turbines, enabling a more precise evaluation of wind farm energy output.
Resolving the detailed boundary layers around turbine blades in CFD simulations requires significant computational resources, posing a major challenge for large-scale wind farm simulations. 
To address this issue, simplified actuator-based models, such as the actuator line model (ALM) and actuator disk model (ADM), have been introduced to reduce computational complexity.
Both models have been shown to accurately capture wake effects in wind farm simulations \citep{martinez2015large,stevens2018comparison,tu2023aerodynamic}, with ADM offering higher computational efficiency compared to ALM. 
In an effort to provide higher fidelity results in WFLO, ADM has recently been used in a number of studies \citep{antonini2020optimal,cruz2020wind,wang2024optimization}. 
Nevertheless, despite ADM's enhanced efficiency, the computational cost of CFD-based WFLO remains significantly higher than that of conventional wake-model-based methods, which continues to limit the widespread adoption of high-fidelity approaches.

Surrogate-based optimization (SBO) is an effective strategy to address the computational challenges of using the CFD-based methods in optimization problems.
The core concept of SBO is to leverage limited high-fidelity samples to construct a surrogate model of the complex objective function, thereby accelerating the optimization process.
The Kriging model, also known as Gaussian process regression, is the most widely used due to its ability to capture the nonlinear relations with fewer samples \citep{gaussian2004}.
Surrogates built upon initial samples are usually not accurate enough to yield reliable optimized results.
As a result, SBO is typically augmented with a Bayesian approach by infilling additional samples to iteratively update the surrogate model until the optimization process converges to the optimal solution.
The key challenge of designing the infilling criterion lies in balancing the trade-off between exploitation and exploration \citep{frazier2018tutorial}. 
Exploitation focuses on refining the search in regions of the design space where high-performing solutions have already been identified, while exploration seeks to query less-explored or uncertain regions to discover potentially better solutions.
Bayesian optimization frameworks often employ acquisition functions (AF) to manage the competing goals of exploitation and exploration. 

The AFs generally fall into two types: improvement-based and information-based strategies.
The improvement-based acquisition functions prioritize exploiting areas that are likely to yield a better function value than the best observed so far.
These functions guide the selection of the next sample by quantifying the potential benefit of sampling at different locations, ensuring a more systematic search of the design space.
A straightforward approach is to update the surrogate model by selecting the next sample based on the maximum surrogate prediction (MSP). 
Other typical examples include Probability of Improvement \citep{kushner1964new} , which assesses the likelihood of surpassing the current best value, and Expected Improvement \citep{mockus1974bayesian}, which considers the magnitude of the potential improvement. 
The improvement-based AFs tend to be computationally efficient and are well-suited for problems where a good starting point is known. 
However, they tend to favor exploitation over exploration, which can lead to convergence on local optima, especially in complex or highly multimodal search spaces. \citep{shahriari2015taking,qin2017improving,wang2024optimization}.

Information-based acquisition functions, on the other hand, are designed to maximize the knowledge gained about the objective function, emphasizing exploration and aiming to reduce uncertainty in the model's predictions. 
They prioritize areas that, once sampled, would reveal the most about the function's landscape, even if the immediate potential improvement isn't maximal.
Entropy Search \citep{hennig2012entropy} is a pioneering information-based AF, whose goal is to maximize the information \citep{gallager1968information} gained about the objective function by reducing uncertainty (specifically, by minimizing the entropy) of the posterior distribution over the global minimum's location.
Predictive Entropy Search \citep{hernandez2016predictive} achieves similar levels of exploration while being significantly more computationally feasible by focusing on reducing the uncertainty in the model's predictions, rather than the uncertainty of the global optimum's location directly.
Max-value Entropy Search \citep{wang2017max} is a variant of Entropy Search that aims to maximize the information gain about the function's optimal value, rather than the locations of the optimum.
It offers more robustness and potentially faster computation compared to standard Entropy Search when the focus is primarily on optimizing the best value and the landscape has multiple local optima.
Despite the strong exploration capabilities, these information-based acquisition functions typically require a greater number of sample evaluations to converge compared to their improvement-based counterparts. 
Moreover, their computational cost often escalates considerably with increasing dimensionality, rendering them less efficient for high-dimensional optimization problems, as is the case in wind farm layout design \citep{thomas2023comparison}.

In summary, while a variety of AFs have been proposed for Bayesian optimization (BO), there remains ongoing debate regarding which one delivers the best performance \citep{shahriari2015taking,qin2017improving,rehbach2020expected}.
Although the improvement-based methods often offer a good balance between exploration and exploitation, they can struggle with complex, multimodal functions, potentially converging to local optima. 
Conversely, the information-based methods are capable of more robust global exploration, but they often come with a higher computational cost, making them impractical for resource-constrained scenarios or high-dimensional problems. 
The performance is often highly sensitive to the hyperparameters within each method, such as the exploration-exploitation trade-off parameters and the Gaussian kernel parameters. 
The effectiveness of these AFs is also heavily influenced by the characteristics of the objective function being optimized, such as its dimensionality, modality (number of local optima), smoothness, etc. 
These factors emphasize the importance of carefully selecting acquisition functions tailored to each specific optimization task and underscore the necessity for adaptive strategies capable of dynamically adjusting to the problem's landscape.

This paper introduces an adaptive switching strategy for acquisition functions (AFs) within a Bayesian optimization framework, specifically designed to facilitate the efficient design of wind farm layouts using high-fidelity methods, which have not been widely used due to the high computational cost.
Among the few who have conducted such analysis, \citet{cruz2020wind} and \citet{antonini2020optimal} applied genetic algorithm (GA) and adjoint method to directly optimize wind farm layout using CFD-based methods instead of wake models.
Recognizing the immense computational demand of the direct optimization, \citet{bempedelis2024data} applied high-fidelity SBO with lower confidence bound (LCB) as the acquisition function.
In our previous study \citep{wang2024optimization}, SBO with MSP as the acquisition function was used.
While this approach is able to significantly reduce the number of CFD simulations, the obtained AEP is often smaller than the direct GA.
This is primarily because MSP exclusively focuses on exploitation rather than sufficient exploration. 
This inherent limitation motivates the development of an adaptive AF switching strategy.

A key contribution of this study is the synergistic integration of two distinct AFs: Maximum Surrogate Prediction (MSP) and Max-value Entropy Search (MES). 
Initially, MSP, with its inherent emphasis on exploitation, is employed to rapidly locate a set of promising optimal solutions. 
Subsequently, leveraging the knowledge acquired during the MSP phase, MES assumes control to ensure that the optimization process effectively escapes local optima, thereby enhancing the global optimization capabilities of the Bayesian framework.
Another contribution of this paper is that an improved sampling method is proposed to generate initial samples for constructing the Kriging model for wind farms with irregular boundaries.
The methods presented in this research achieve a substantial reduction in AEP evaluations within surrogate-based optimization for wind farm layout design, allowing for the integration of high-fidelity CFD simulations into WFLO investigations.

\section{Bayesian optimization framework for wind farm layout design}
\label{sec:BO}
In this section, the details about the Bayesian optimization framework for wind farm layout optimization are presented.
The construction of Kriging for Bayesian optimization is introduced in Section \ref{sec:kriging}.
Next, in Section \ref{sec:lhs}, a Latin-Hypercube-Sampling-based (LHS-based) constrained sampling method is proposed to generate the initial samples for the optimization problems with irregular design space.
We further outline the proposed adaptive switch strategy for acquisition functions for the \emph{exploration-exploitation} balance in Section \ref{sec:af}.
The penalty function approach is then introduced to handle the constraint in Section \ref{sec:penalty}. 
The termination criterion is then illustrated, and the Bayesian optimization framework is formulated in Section \ref{sec:boframework}.
Finally, we perform optimization on two analytical functions to evaluate the performance of the BO framework based on the adaptive switch strategy for AFs in Section \ref{sec:experiments}.

\subsection{Kriging model}
\label{sec:kriging}
The Kriging model \citep{krige1951statistical,Rasmussen2006Gaussian}, also known as Gaussian Process (GP), is a non-parametric probabilistic model.
It represents an approximation to the design function $f(\boldsymbol{x})$ as a combination of a trend function and a stochastic process:
$\epsilon(\boldsymbol{x})$ as
\begin{equation}
    \hat{f}(\boldsymbol{x}) = \sum_{i=1}^k\beta\textit{g}_i(\boldsymbol{x}) + \epsilon(\boldsymbol{x}),
    \label{equ:KG}
\end{equation}
where $\hat{f}(\boldsymbol{x})$ is the Kriging approximation to the design function, $\boldsymbol{x}$ is the design variable.
The trend functions, represented by $\Sigma_{i=1}^k\beta\textit{g}_i(\boldsymbol{x})$, are assumed to be constant in this study.
The error term, $\epsilon(\boldsymbol{x})$, is modeled as a Gaussian process, characterized by zero uncertainty at the training points \citep{adams2020dakota}.
The covariance structure of this Gaussian process is defined as
\begin{equation}
    \mathrm{Cov}[\epsilon(\boldsymbol{x}^{(i)}),\epsilon(\boldsymbol{x}^{(j)})] = \sigma^{2}R(\boldsymbol{x}^{(i)},\boldsymbol{x}^{(j)}),
\end{equation}
where $\textit{R}$ is the correlation function.
In this study, the power exponential correlation is used:
\begin{equation}
     R(\boldsymbol{x}^{(i)},\boldsymbol{x}^{(j)}) = \mathrm{exp}(-\sum_{k=1}^m\theta_k|{x}^{(i)}_{k}-{x}^{(j)}_{k}|^{p}), \ \ \ i,j\ =\ 1,2,...,n. 
\end{equation}
Training a Kriging model involves determining the kernel hyperparameters, $\theta_k$ and $p$, within the correlation function; further details can be found in \citet{KrigingTrain}. For this study, we use the SMT 2.0 Python package \citep{saves2024smt} to perform Kriging model training.

\subsection{Constrained Latin hypercube sampling method}
\label{sec:lhs}
Initial samples are required to build the Kriging model.
The Latin hypercube sampling (LHS) is the most widely used random sampling method for generating the initial samples \citep{mckay2000comparison}.  
LHS aims to ensure that the sample points are evenly distributed across the rectangle domain of each input variable, thereby improving the representation of the design space.
For a function of $k$ random variables, LHS is performed to obtain $N$ samples in the following manner:

\begin{enumerate}
    \item Divide the rectangle domain: The domain of each variable $X_i$, $i = 1,2,...,k$ is divided into $N$ non-overlapping intervals of equal probability.
    \item Sample each interval: One sample is randomly selected from each interval for each variable.
    \item Combine samples: The sampled values for each variable are then randomly paired to create $N$ sample points. This ensures that each interval for each variable is represented in the final sample set.
\end{enumerate}

Standard LHS method assumes a rectangular sampling domain defined by the bounds of each variable. 
However, in WFLO problems, the wind farm boundaries can be irregular, making it difficult to apply standard LHS. 
To address this, a constrained LHS method that accommodates irregular wind farm boundaries is proposed.
This method first generates a large number of samples using standard LHS, encompassing the entire wind farm area. 
Then, it filters these samples, retaining only those that fall within the irregular boundaries. 
If the number of feasible samples is insufficient, an expanded LHS technique \citep{arenzana2021multi,saves2024smt} is applied to add new samples while preserving the space-filling properties of LHS, ensuring a representative set of initial samples for the Kriging model.
A pseudo code for the constrained LHS method is provided in Algorithm \ref{algo:lhs}.

\begin{algorithm}[H]
    \caption{Constrained LHS method}
    \label{algo:lhs}
    \KwIn{The rectangle domain defined by variables' bounds $\boldsymbol{D_0}$, number of the required samples $n$, constrain $F$ }
    \KwOut{A group of samples within the irregular design space, $\boldsymbol{X}$=$\{\boldsymbol{X}_1$, $\boldsymbol{X}_2$, ..., $\boldsymbol{X}_N\}$}
    Define the irregular design space $\boldsymbol{D}$ with $\boldsymbol{D_0}$ and $F$\;
    Initialize $m$ samples within  $\boldsymbol{D}_0$ using the standard LHS ($m \gg n$)\;
    Select initialized samples within $\boldsymbol{D}$ and add them to $\boldsymbol{X}$\;
    $n_0 \leftarrow $ number of samples within $\boldsymbol{X}$\;
    \While{$n_0$ < $n$}{
    	Expand $k$ samples apart from $\boldsymbol{X}$ using the expanded LHS ($k=n-n_0$)\;
	Select expanded samples within $\boldsymbol{D}$ and add them to $\boldsymbol{X}$\;
	$n_0 \leftarrow $ number of samples within $\boldsymbol{X}$\;
	}
\end{algorithm}

\subsection{The adaptive switch strategy for acquisition functions}
\label{sec:af}
A surrogate model built on the initial dataset is usually insufficient for accurate optimization. 
More samples are required to refine the surrogate and enable a better search for the optimum.
In order to select the next sampling point, various acquisition functions (AF) have been developed, and candidate sample $\boldsymbol{x}_{n+1}$ for the next evaluation of $f$ is selected in iteration $n$ by optimizing these acquisition functions as
\begin{equation}
    \boldsymbol{x}_{n+1} = \mathrm{argmax}(AF_n(\boldsymbol{x})).
\end{equation}

\subsubsection{Maximum Surrogate Prediction}
\label{sec:MSP}
The Maximum Surrogate Prediction (MSP) is the most basic improvement-based acquisition function (AF) which feeds the optimized results back to the dataset sequentially. 
MSP is formally defined as:
\begin{equation}
	\mathrm{AF}_{\mathrm{MSP},n}(\boldsymbol{x}) = \hat{f}_n(\boldsymbol{x}),
\end{equation}
where $\hat{f}_n(\boldsymbol{x})$ denotes the Kriging emulator built in the $n$-th iteration. 
This acquisition function fully exploits available data without taking the surrogate's uncertainty into consideration.
While this approach can quickly locate promising areas, it often suffers from premature convergence, as the search may be stuck in local optima.

\subsubsection{Max-value Entropy Search}
\label{sec:MES}
The Max-value Entropy Search (MES) \citep{wang2017max} acquisition function quantifies the expected information gain about the maximum value of the objective function. 
It calculates the reduction in entropy of the distribution of the maximum value after incorporating a new observation.
The acquisition function of the MES is formulated as follows:
\begin{equation}
    \mathrm{AF}_{\mathrm{MES},n}(\boldsymbol{x}) = H(p(y^*|D_n)) - H(p(y^*|D_n,\boldsymbol{x},y)),
\end{equation}
where $H(p(\boldsymbol{x}))=-\int p(\boldsymbol{x})\mathrm{log}p(\boldsymbol{x})\mathrm{d}\boldsymbol{x}$ is the information entropy \citep{gray2011entropy} which quantifies the level of uncertainty or information content associated with the probability distribution $p(\boldsymbol{x})$.
$p(y^*|D_n)$ is the probability measure on $y^*$, which is the $max$ of the objective function inferred by the Kriging model built upon samples $D_n$.
After incorporating sample $\boldsymbol{x}$, the probability measure on $y^*$ is changed to $p(y^*|D_n,\boldsymbol{x},y)$.
However, it is hard to directly calculate the entropy of $p(y^*|D_n,\boldsymbol{x},y)$ due to the existence of random variables $y$.
Its expectation is utilized here as follows:
\begin{equation}
    \mathbb{E}_{y}[H(p(y^*|D_n,\boldsymbol{x},y))] = \int p(y|D_n,\boldsymbol{x})H(p(y^*|D_n,\boldsymbol{x},y))\mathrm{d}y, 
\end{equation}
The $\mathrm{AF}_{\mathrm{MES},n}(\boldsymbol{x})$ then transforms into:
\begin{subequations}
\label{equ:mes_af}
\begin{align}
    \mathrm{AF}_{\mathrm{MES},n}(\boldsymbol{x}) &= H(p(y^*|D_n))-\mathbb{E}_{y}[H(p(y^*|D_n,\boldsymbol{x},y))] \\
    & = H(p(y|D_n,\boldsymbol{x}))- \mathbb{E}_{y^*}[H(p(y|D_n,\boldsymbol{x},y^*))].
\end{align}
\end{subequations}
The first term in equation \ref{equ:mes_af} (b) is the entropy of the Gaussian distribution $p(y|D_n,\boldsymbol{x})$ which is computed in closed form as shown below:
\begin{equation}
    H(p(y|D_n,\boldsymbol{x})) = \frac{1+\mathrm{ln}(2\pi)}{2} + \mathrm{ln}(\sigma(\boldsymbol{x})),
\end{equation}
where $\sigma(\boldsymbol{x})$ is the Kriging's predicted variances at sample $\boldsymbol{x}$.
For the second term, the Monte Carlo sampling approach is utilized to approximate as:
\begin{equation}
    \mathbb{E}_{y^*}[H(p(y|D_n,\boldsymbol{x},y^*))] \approx \frac{1}{K}\sum^K_{i=1}[H(p(y|D_n,\boldsymbol{x},y^*_i))],
\end{equation}
where $K$ is the number of sampled optimums $y^*$.
The $p(y|D_n,\boldsymbol{x},y^*_i)$ is a truncated Gaussian with upper bound $y^*_i$, and the differential entropy for it is given as:
\begin{equation}
   H(p(y|D_n,\boldsymbol{x},y^*_i)) \approx  \frac{1+\mathrm{ln}(2\pi)}{2}+\mathrm{ln}(\sigma(\boldsymbol{x}))+\mathrm{ln}\Phi(\gamma_{y^*}(\boldsymbol{x}))-\frac{\gamma_{y^*}(\boldsymbol{x})\phi(\gamma_{y^*}(\boldsymbol{x}))}{2\Phi(\gamma_{y^*,i}(\boldsymbol{x}))},
\end{equation}
where $\gamma_{y^*,i}(\boldsymbol{x}) = (y^*_i-\mu(\boldsymbol{x}))/\sigma(\boldsymbol{x})$ is the standardized Gaussian distribution, $\phi (\bullet)$ and $\Phi (\bullet)$ represent the probability and cumulative density function of the Gaussian distribution respectively.
Then, the MES AF can be approximated as follows:
\begin{subequations}
\label{equ:mes_af_approx}
\begin{align}
    \mathrm{AF}_{\mathrm{MES},n}(\boldsymbol{x}) &\approx \left[ \frac{1+\mathrm{ln}(2\pi)}{2} + \mathrm{ln}(\sigma(\boldsymbol{x}))\right] - \frac{1}{K}\sum^K_{i=1} \left [ \frac{1+\mathrm{ln}(2\pi)}{2}+\mathrm{ln}(\sigma(\boldsymbol{x}))+\mathrm{ln}\Phi(\gamma_{y^*,i}(\boldsymbol{x}))-\frac{\gamma_{y^*,i}(\boldsymbol{x})\phi(\gamma_{y^*,i}(\boldsymbol{x}))}{2\Phi(\gamma_{y^*,i}(\boldsymbol{x}))}\right] \\
    & = \frac{1}{K}\sum^K_{i=1}\left[ \frac{1+\mathrm{ln}(2\pi)}{2} + \mathrm{ln}(\sigma(\boldsymbol{x}))\right] - \frac{1}{K}\sum^K_{i=1} \left [ \frac{1+\mathrm{ln}(2\pi)}{2}+\mathrm{ln}(\sigma(\boldsymbol{x}))+\mathrm{ln}\Phi(\gamma_{y^*,i}(\boldsymbol{x}))-\frac{\gamma_{y^*,i}(\boldsymbol{x})\phi(\gamma_{y^*,i}(\boldsymbol{x}))}{2\Phi(\gamma_{y^*,i}(\boldsymbol{x}))}\right]\\
    & = \frac{1}{K} \sum^K_{i=1} \left[\frac{\gamma_{y^*,i}(\boldsymbol{x})\phi(\gamma_{y^*,i}(\boldsymbol{x}))}{2\Phi(\gamma_{y^*,i}(\boldsymbol{x}))}-\mathrm{ln}(\Phi(\gamma_{y^*,i}(\boldsymbol{x})))\right]
\end{align}
\end{subequations}

\subsubsection{The adaptive switch-AF strategy}
\begin{figure}
    \centering
    \includegraphics[width=1.0 \textwidth]{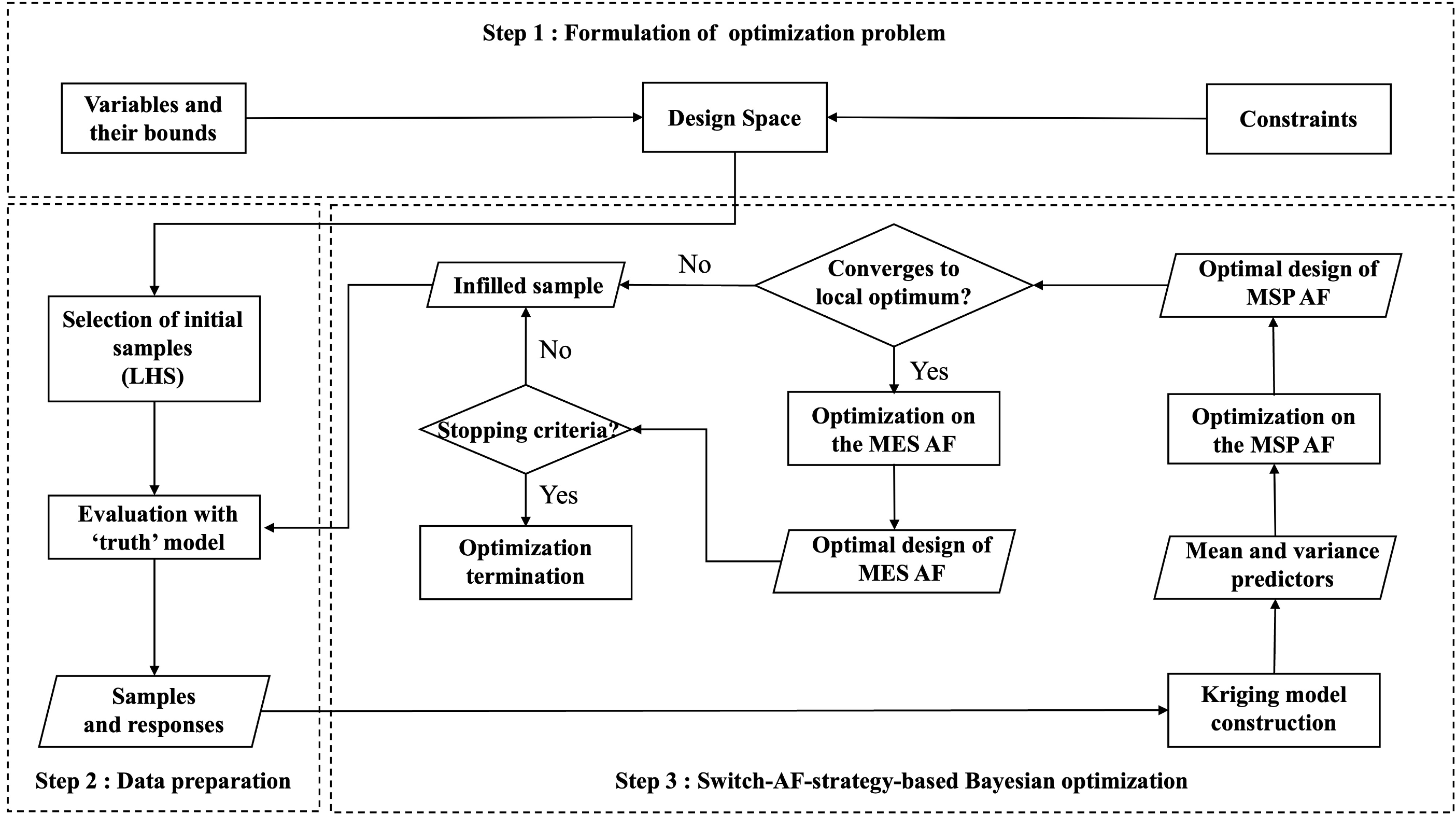}
    \caption{The workflow of Bayesian optimization using the proposed adaptive switching strategy between MSP and MES AFs.}
    \label{fig:bo}
\end{figure}
The MES equation \ref{equ:mes_af_approx} relies heavily on Kriging model to infer the value of optima, $y^*$.
However, it is typically hard to find appropriate initial samples which provide meaningful prior knowledge of the landscape of the design space.
Searching for the infilling samples using only MES may result in the waste of computational resources in the initial phase of optimization.
It is reminded that MES does not exploit the inferred optimal location $\boldsymbol{x}^*$ in its formulation and sampling strategy, as discussed in Section \ref{sec:MES}.  
On the other hand, MSP actively incorporates the inferred optimal location $\boldsymbol{x}^*$ as a new sample; this emphasis on exploitation allows MSP to locate high-performing solutions efficiently.
For this reason, an adaptive switch strategy between MSP and MES is proposed in this study to leverage their respective strengths. 
The MSP is first used to find samples near local optima in each iteration of the BO.
Once the inferred optimum $x^*$ found by MSP AF converges, the MES AF is employed to refine the search for global optima.
The switch criterion from MSP to MES is defined as
\begin{equation}
\label{equ:msp_convergence}
    \sum^d_{j=1} |x^*_j - X_{i,j}| \leq \epsilon,\ \ \ \ i=1,2,...,M,
\end{equation}
where $x^*_j$ is the $j$-th variable in the optimum $\boldsymbol{x}^*$, $X_{i}$ is the $i$-th training sample and $\epsilon$ is a pre-defined value (0 for discrete optimization problem).
The adaptive switching mechanism ensures that the optimization maintains a balance between efficiency and global exploration capability.
The workflow of Bayesian optimization utilizing the adaptive switching strategy for AFs is depicted in Figure \ref{fig:bo}.

For the optimization algorithm, the differential evolution (DE) algorithm in the Python library \verb|scipy| package \citep{2020SciPy-NMeth} is implemented due to its computational efficiency.
The optimization on MSP AF is implemented in parallel, in order to overcome the influence of initial settings and obtain a group of $y^*$ to feed into MES AF.
Furthermore, in order to enhance the optimization performance and efficiency of the BO algorithm, a population initialization and mutation strategy based on existing samples is introduced in this study for DE as illustrated below:
\begin{enumerate}
    \item Randomly generate a set of samples $\boldsymbol{X}_\mathrm{random}$ using the constrained sampling method.
    \item Perform random mutation on variables of the best-performed sample to obtain a new group of samples $\boldsymbol{X}_\mathrm{mutate}$.
     \item Combine the two sets of samples above and provide them for DE as the initial population.
\end{enumerate}

\subsection{Penalty function}
\label{sec:penalty}
A penalty function approach is employed to handle the complex constraints during the optimization task. 
This is achieved by modifying the optimization objective when the constraint inequality is violated for a given input $\boldsymbol{X}$. 
The transformed objective function is expressed as:
\begin{equation}
    \mathrm{max}\ \ \ \  \lambda(\boldsymbol{X})\mathrm{AEP}(\boldsymbol{X}),
\end{equation}
where the $\lambda(\boldsymbol{X})<1$ is the penalty factor that  effectively reduces the objective value when constraints are violated.

\subsection{Termination criterion and formulation of the Bayesian optimization framework}
\label{sec:boframework}
Three termination criteria are defined in the current BO framework. 
The first one is on the computational budget:
\begin{equation}
	N \leq N_\mathrm{max},
\end{equation}
where $N$ is the number of function calls.

The second one is to monitor the diversity of the training samples. Similar to the approach used in genetic algorithms, the Bayesian optimization also defines a population, $\boldsymbol{X}_{\mathrm{pop}} = \{\boldsymbol{x}_1,\boldsymbol{x}_2,...,\boldsymbol{x}_n\}$. These $n$ samples are the top-performing ones within the training samples. The deviation of their objective values, $var$, should be less than a pre-defined value $\epsilon_1$ as:
\begin{equation}
	var \leq \epsilon_1.
\end{equation}

The third one is specific for the MES acquisition function, i.e., 
\begin{equation}
	\mathrm{max}(\mathrm{AF}_{\mathrm{MES}}) \leq \epsilon_2,
\end{equation}
where $\epsilon_2$ is a user-defined limit.
Once any one of these criteria above is reached, the optimization is stopped.

With the Kriging model, constrained LHS method, adaptive switch strategy for AF, penalty function approach, and the termination conditions illustrated above, the BO framework for the WFLO problem is formulated. 
The procedure is described in Algorithm \ref{algo:bo}. 
More detailed information on the implementation of the BO framework is provided in the github repositories\footnote{\url{https://github.com/sjtuwzf/WFLBO}}.

\begin{algorithm}[H]
    \caption{The Bayesian optimization framework}
    \label{algo:bo}
    \KwIn{The space defined by variables' bounds $\boldsymbol{D_0}$, constrain $\boldsymbol{F}$,  computationally expensive objective function $f$. }
    \KwOut{The optimized results $\boldsymbol{x}_{\mathrm{opt}}$ and its objective function value $y_{\mathrm{opt}}$}
    Define the irregular design space $\boldsymbol{D}$ with $\boldsymbol{D_0}$ and $F$\;
    Initialize training samples  $\boldsymbol{X}_0 = \{\boldsymbol{x}_1,\boldsymbol{x}_2,...\boldsymbol{x}_m\}$ using the constrained sampling method\;
    Evaluate these samples with $f$ and obtain $\boldsymbol{Y}_0 = \{\boldsymbol{y}_1,\boldsymbol{y}_2,...\boldsymbol{y}_m\}$\;
    \For{$i = 0,1,2,3,...,N_{\mathrm{max}}$}{
   	Select population $\boldsymbol{X}_{\mathrm{pop}}$ and $\boldsymbol{Y}_{\mathrm{pop}}$\;
	Compute the standard deviation $var$ of $\boldsymbol{Y}_{\mathrm{pop}}$\;
	\If{$var \geq \epsilon_1$}{
		Build Kriging emulator $\hat{f}$ with $\boldsymbol{X}_i$ and $\boldsymbol{Y}_i$\;
		$\boldsymbol{x}_i\leftarrow \mathrm{argmax}_{x\in \boldsymbol{D}} \hat{f}(\boldsymbol{x})$\;
		\If{$\boldsymbol{x}_i\ \mathrm{converges}$}{
			$\boldsymbol{x}_i \leftarrow \mathrm{argmax}_{\boldsymbol{x}\in \boldsymbol{D}} \mathrm{AF}_{\mathrm{MES},i}(\boldsymbol{x})$\;
			\If{$\mathrm{AF}_{\mathrm{MES},i}(\boldsymbol{x}_i) < \epsilon_2$}{
				Terminate optimization\;
			}
		}
		$y_i \leftarrow f(\boldsymbol{x}_i)$\;
		$\boldsymbol{X}_i \leftarrow  \boldsymbol{X}_{i-1}\cup \boldsymbol{x}_i$\;
		$\boldsymbol{Y}_i \leftarrow \boldsymbol{Y}_{i-1}\cup y_i$\;
	}
	\Else{
		Terminate optimization\;
	}
    }
    $(\boldsymbol{x}_{\mathrm{opt}},y_{\mathrm{opt}}) \leftarrow (\mathrm{argmax}\boldsymbol\ {Y},\mathrm{max}\boldsymbol\ {Y})$
\end{algorithm}

\subsection{Benchmark using the Ackley Function}
\label{sec:experiments}
\begin{figure}
    \centering
    \includegraphics[width=0.7 \textwidth]{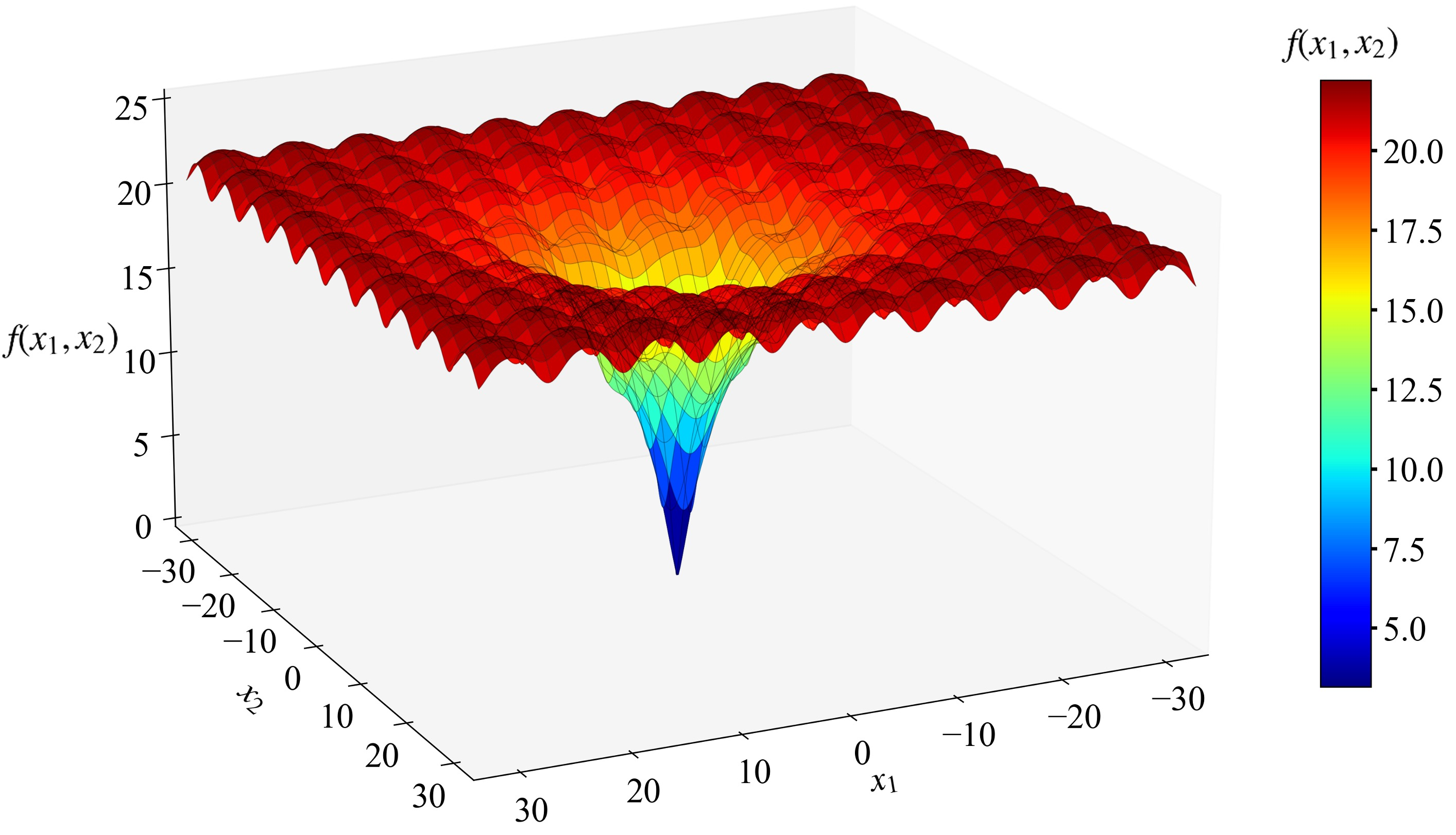}
    \caption{The response surface of the 2 dimensional Ackley function}
    \label{fig:2dAckley}
\end{figure}

In this section, the performance of the adaptive switch strategy between MSP and MES acquisition functions is tested on two challenging optimization test functions: the 4-dimensional (4D) and the 10-dimensional (10D) Ackley functions.
The Ackley function is expressed as
\begin{equation}
	f(\boldsymbol{x}) = -a\mathrm{exp}\left(-b\sqrt{\frac{1}{d}\sum^d_{i=1}x^2_i}\right)-\mathrm{exp}\left(\frac{1}{d}\sum^d_{i=1}\mathrm{cos}(cx_i)\right)+a+\mathrm{exp}(1),
\end{equation}
where $d$ is the dimension of the problem. The optimization goal is to find the minimum of the functions.
The Ackley function is widely used for testing optimization algorithms due to the existence of many local optima, as illustrated in Figure \ref{fig:2dAckley}.
The constants are set as $a = 20$, $b = 0.2$, and $c=2\pi$ as suggested by \citet{motiian2011}, and the function is evaluated on the hypercube $x_i \in [-32.768,32.768]$, for all $i=1,...,d$.
The global minimum is given as $f(\boldsymbol{x}^*) = 0$, at $\boldsymbol{x}^* = (0,0,...,0)$.
The standard LHS method is used to produce 10 and 25 initial samples for the two functions.

\begin{figure}
    \centering
    \includegraphics[width=0.8\textwidth]{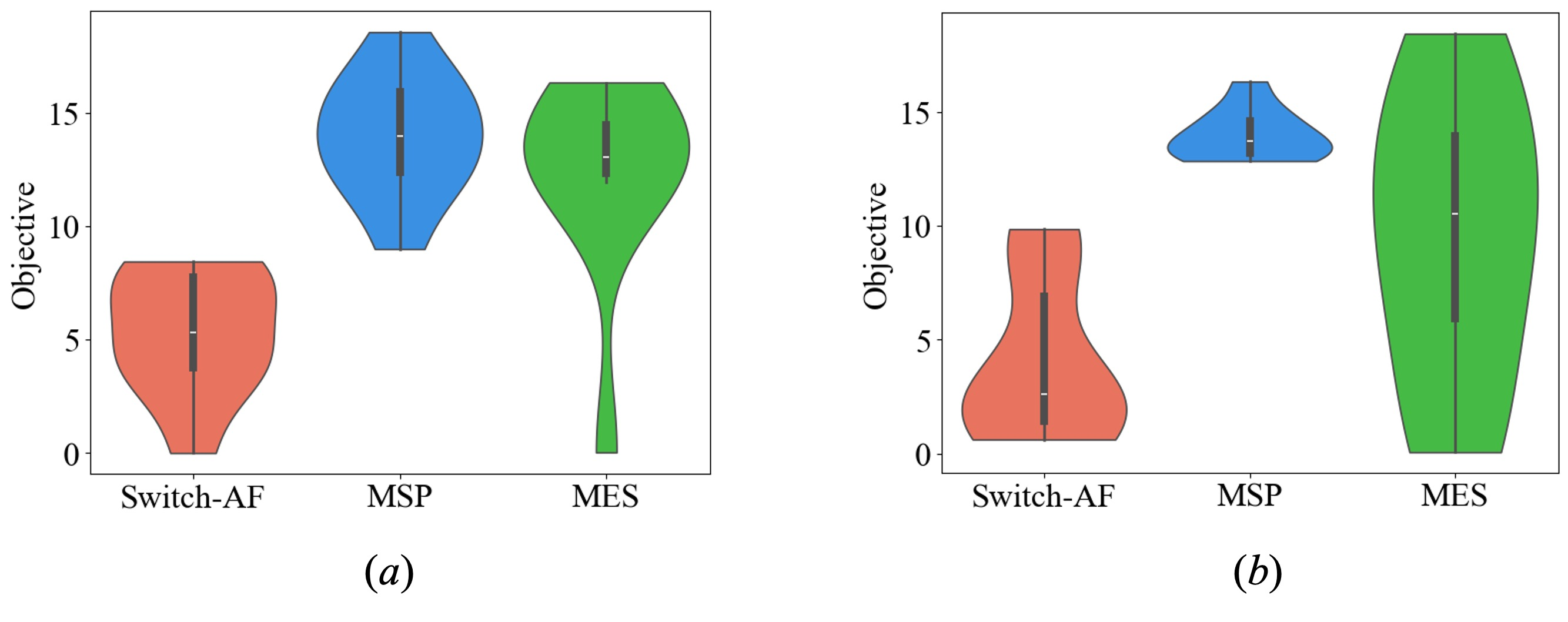}
    \caption{The violin plots of the optimized results obtained by using switch-AF, MSP and MES for ($a$) the 4D Ackley function and ($b$) the 10D Ackley function.}
    \label{fig:ackley_violin}
\end{figure}

\begin{figure}
    \centering
    \includegraphics[width=1.0\linewidth]{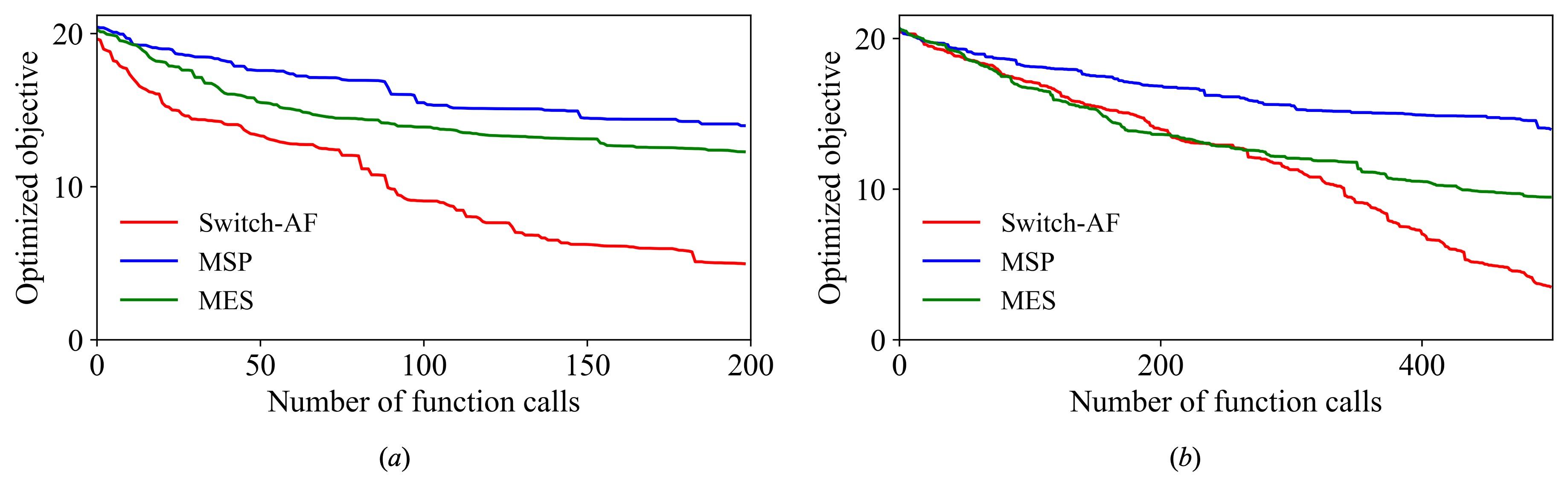}
    \caption{The development of the optimized results obtained by using switch-AF, MSP and MES for ($a$) the 4D Ackley function and ($b$) the 10D Ackley function.}
    \label{fig:ackley_iterations}
\end{figure}

The optimizations are performed using three acquisition functions, i.e., MSP, MES and the proposed adaptive switch between the two, with a budget of 200 and 500 function evaluations for the 4D and 10D problems, respectively. 
Each optimization is repeated 10 times with randomly generated training samples to average out the initial condition effects.
The optimization results, visualized through violin and iteration evolution as shown in Figures \ref{fig:ackley_violin} and \ref{fig:ackley_iterations}, demonstrate the superiority of the proposed adaptive switch acquisition function (switch-AF) compared to both MSP and MES. 
While MES exhibits the potential for outstanding performance in isolated instances, even approaching the global optimum (0), its overall efficacy was generally lower than the switch-AF. 
This highlights a key limitation of MES under limited budgets: its reliance on random sampling may not efficiently locate optimal regions. 
On the other hand, MSP also exhibits less satisfying results due to its susceptibility to converging to local optima prematurely.
The strength of the switch-AF lies in its strategic combination of these approaches: it initially leverages MSP's fast convergence to gain improved prior knowledge of the optimum, subsequently transitioning to MES for broader global exploration once MSP's initial gains plateau. 
It is also noted that the cost of constructing the Kriging model increases with the number of samples used during the optimization process, making it impractical to conduct a large number of iterations. 
Nonetheless, the superior optimization performance of the proposed switch-AF makes it a promising technique for addressing computationally expensive black-box optimization problems.

\section{Wind farm layout optimization}
\label{sec:WFLO}

\subsection{Formulation of the WFLO problem}
\label{sec:problem}
In the WFLO problem selected for this study, the objective is to maximize the AEP of a wind farm by finding the optimal positions of each wind turbine  $[\boldsymbol{x},\boldsymbol{y}]^T$, with $\boldsymbol{x} = (x_1,x_2,...,x_N)$ and $\boldsymbol{y} = (y_1,y_2,...,y_N)$, where $N$ is the number of turbines within the wind farm, and $(x_i,y_i)$ represents the coordinate of the $i$-th turbine. 
The AEP of the wind farm is calculated as the weighted sum of each turbine's power under different wind conditions (i.e., wind speed and wind direction):
\begin{equation}
    \mathrm{AEP}=\frac{365\times3+366}{4}\times24\sum_{i=1}^{n}\sum_{j=1}^{m}P_{j}(x_{i},y_{i})f_{j},
\end{equation}
where  $P_{j}(x_{i},y_{i})$ is the turbine's mechanical power output under wind state $j$, with $f_{j}$ representing its frequency of occurrence.

Wind turbines are only allowed to be installed within the pre-specified wind farm boundary for the WFLO problem. 
The constraints for the WFLO problem are therefore expressed as follows:$(x_i,y_i) \in \boldsymbol{D}$, where $\boldsymbol{D}$ denotes the design space defined by the wind farm's boundary.

\subsection{Gaussian wake model}
\label{sec:gch}
Although the aim of this study is to enable high-fidelity WFLO, the low-fidelity Gaussian wake model is used here to demonstrate the feasibility of the proposed framework.
The Gaussian wake model assumes a normal distribution of velocity profiles in the wake region \citep{king2021control}.
The streamwise velocity behind a turbine, $u_G$ is computed as follows:
\begin{subequations}
\begin{align}
	& \frac{u_G}{U_\infty} = 1 - C\mathrm{exp}\left(\frac{-(y-y_0)^2}{2\sigma^2_ye^{-(z-z_h)^2/2\sigma_z^2}}\right),\\
	& C = 1-\sqrt{1-\frac{(\sigma_{y0}\sigma_{z0})C_T}{\sigma_y\sigma_z}},
\end{align}
\end{subequations}
in which $U_\infty$ is the freestream velocity, $C$ is the wake center's velocity deficit, $\sigma_y$ and $\sigma_z$ represent the wake width in $y$ and $z$ direction, $y_0$ is the turbine's spanwise position and $z_h$ is the hub height. 
$\sigma_y$ and $\sigma_z$ are computed as
\begin{subequations}
\begin{align}
	& \frac{\sigma_z}{D} = k_z\frac{(x-x_0)}{D}+\frac{\sigma_{z0}}{D},\\
	& \frac{\sigma_y}{D} = k_y\frac{(x-x_0)}{D}+\frac{\sigma_{y0}}{D},
\end{align}
\end{subequations}
where $\frac{\sigma_{z0}}{D}=\frac{\sigma_{y0}}{D} = \frac{1}{2}\sqrt{\frac{u_R}{U_\infty+u_0}}$, $D$ is the rotor diameter, $u_R$ is the velocity at the rotor, $u_0$ is the velocity of the far wake, $k_y$ and $k_z$ define the wake expansion in the lateral and vertical directions.

Once the velocity at each turbine is computed, the power coefficient $C_P$ is obtained by looking up the table provided by the wind turbine manufacturers. 
With the power coefficient $C_P$, wind turbines' mechanical power is calculated as follows:
\begin{equation}
P = \frac{1}{2}\rho AU^3_\infty C_P
\end{equation}
where $P$ is the produced mechanical power, $\rho$ is the air density.
In this study, the Gaussian wake model is implemented in the open-source Python library \verb|Floris| \citep{Floris}, which is developed by the National Renewable Energy Laboratory (NREL) for evaluating the wake losses among wind farms.

\subsection{CFD-based wind farm energy output prediction}
\label{sec:cfd}
In this study, the Reynolds-averaged Navier-Stokes (RANS) formulation is employed to model the wake interactions within a wind farm. Similar to other  CFD-based WFLO studies \citep{cruz2020wind,antonini2020optimal}, the rotors are simplified as actuator disks. The governing equations are illustrated as follows:
\begin{subequations}
    \label{equ:rans}
    \begin{align}
          \frac{\partial \overline{u}_i}{\partial x_i} & = 0, \\
         \rho \frac{\partial \overline{u}_i}{\partial t} + \rho\frac{\partial(\overline{u}_i\overline{u}_j)}{\partial x_j} & = -\frac{\partial \overline{p}}{\partial x_i}+\frac{\partial}{\partial x_j}\left(\mu\left(\frac{\partial \overline{u}_i}{\partial x_j}+\frac{\partial \overline{u}_j}{\partial x_i}\right)-\rho\overline{u^{'}_{j} u^{'}_{i}}\right) + f_i,
    \end{align}
\end{subequations}
where $x_{i}$ is the Cartesian space coordinate, $\overline{u}_{i}$ and $\overline{p}$ are the temporally averaged velocity and pressure, respectively.
$\rho$ and $\mu$ are the air density and dynamic viscosity, $f_i$ represents the source term from the actuator disk.
By introducing two new variables, i.e., the turbulence kinematic energy $k$ and the dissipation rate $\epsilon$, the $k$-$\epsilon$ turbulence model solves the time-averaged Reynolds stress $\rho\overline{u_{j}^{'}u_{i}^{'}}$.
The standard transport equations for the two new variables are as follows:
\begin{subequations}
    \label{equ:turbulence}
    \begin{align}
        \rho\frac{\partial k}{\partial t}+\rho\frac{\partial (\overline{u}_{i}k)}{\partial x_i} &=\frac{\partial}{\partial x_j}\left(\frac{\mu_t}{\sigma_k}\frac{\partial k}{\partial x_j}\right)+P_k-\rho\epsilon , \\
        \rho\frac{\partial\epsilon}{\partial t}+\rho\frac{\partial(\overline{u}_{i}\epsilon)}{\partial x_i} &=\frac{\partial}{\partial x_j}\left(\frac{\mu_t}{\sigma_\epsilon}\frac{\partial\epsilon}{\partial x_j}\right)+C_{1\epsilon}\frac{\epsilon}{k}P_k - C_{2\epsilon}\frac{\epsilon^2}{k}\rho ,
    \end{align}
\end{subequations}
with 
\begin{subequations}
\label{equ:pkmu}
    \begin{align}
        P_k &= -\rho\overline{u_{j}^{'}u_{i}^{'}}\frac{\partial u_j}{\partial x_i} , \\
        \mu_t &= -\rho C_\mu \frac{k^2}{\epsilon},
    \end{align}
\end{subequations}
where $C_\mu = 0.09$, $\sigma_k=1.0$,  $\sigma_\epsilon=1.3$, $C_{1\epsilon}=1.44$, $C_{2\epsilon}=1.92$ are the five constants in the $k$-$\epsilon$ turbulence model \citep{Launder}.

In the one-dimensional actuator disk model (ADM) \citep{burton2011wind}, thrust and torque of the wind turbine are pre-defined. The volumetric source terms  $f_i$ are calculated by distributing the thrust and torque following Goldstein optimum \citep{goldstein1929vortex} in this study:
\begin{subequations}
\label{equ:distribution}
    \begin{align}
        f_{ix} &= A_{x}r^*\sqrt{1-r^*}, \\
        f_{i\theta} &= A_\theta \frac{r^*\sqrt{1-r^*}}{r^*(1-r^{'}_{h})+r^{'}_{h}},
    \end{align}
\end{subequations}
with
\begin{subequations}
\label{equ:distributionpara}
    \begin{align}
        r^* &= \frac{r^{'}-r^{'}_h}{1-r^{'}_h}, \ \ r^{'} = \frac{r}{R_P}, \ \ r_h^{'} = \frac{R_H}{R_P}, \\
        A_x &= \frac{105}{8}\frac{T}{\pi\Delta(3R_H+4R_P)(R_P-R_H)}, \\
        A_\theta &= \frac{105}{8}\frac{Q}{\pi\Delta R_P(3R_P+4R_H)(R_P-R_H)},
    \end{align}
\end{subequations}
where $f_{ix}$ is axial force, $f_{i\theta}$ is tangential force, $r$ is the distance between the point and disk center, $R_P$ is the external radius of disk, $R_H$ is the internal radius of disk, $T$ is rotor's thrust, $Q$ is rotor's torque, and $\Delta$ is the thickness of disk.

However, it is not straightforward to pre-define the thrust and torque of each rotor, since the theoretical inflow velocity is not known \emph{a priori} due to the wake interactions within the turbine clusters. The thrust and torque curves provide the relationship between inflow velocity and these aerodynamic quantities. Following \citet{richmond2019evaluation}, this study computes the thrust and torque in an iterative manner based on the curves as follows:
 \begin{subequations}
\label{equ:admtheory}
    \begin{align}
        & U_{D} = U_{\infty}(1-a), \\
        & C_{T} = \frac{T}{\frac{1}{2}\rho U^{2}_{\infty}A_{D}}, \\
        & C_{T} = 4a(1-a), \\
        & a = \frac{1}{2}(1-\sqrt{1-C_T}),
    \end{align}
\end{subequations}
where $U_D$ is the velocity averaged over the disk region, $C_T$ is the coefficient of thrust, $A_D$ is the area of the actuator disk, and $a$ is the axial induction factor.

Finally, the power of each turbine is predicted by 
\begin{equation}
\label{equ:power}
    P = \sum_{i=1}^{k}P_{i} = \sum_{i=1}^{k}F_{i}U_{i}
\end{equation}
where $F_i$ and $U_i$ are vectors of volume force and velocity in each individual cell in the actuator disk zone.

The RANS equations are discretized using the finite volume method with second-order numerical schemes in both space and time.  
The simpleFoam solver (OpenFOAM package) \citep{Weller1998OpenFOAM} is modified to incorporate ADM to solve for the steady-state flow over wind turbines, as illustrated above.

\section{Case studies using low-fidelity wake model}
\label{sec:case_studies}
In this section, three wind farm layout optimization cases using low-fidelity wake models are deployed.
The performance of the Bayesian-based WFLO framework is evaluated by comparing it with direct optimization in which the heuristic algorithms including differential evolution (DE), genetic algorithm (GA), and simulated annealing (SA), are directly integrated with Gaussian wake model instead of relying on Kriging model. 
The detailed case descriptions and Bayesian optimization framework settings are illustrated in Section \ref{sec:settings}.
In Section \ref{sec:direct}, the wind farm layouts obtained from three heuristic algorithms and the Bayesian optimization framework are compared. 

\subsection{Case descriptions and optimization
 settings}
\label{sec:settings}

\begin{figure}
    \centering
    \includegraphics[width=1.0\textwidth]{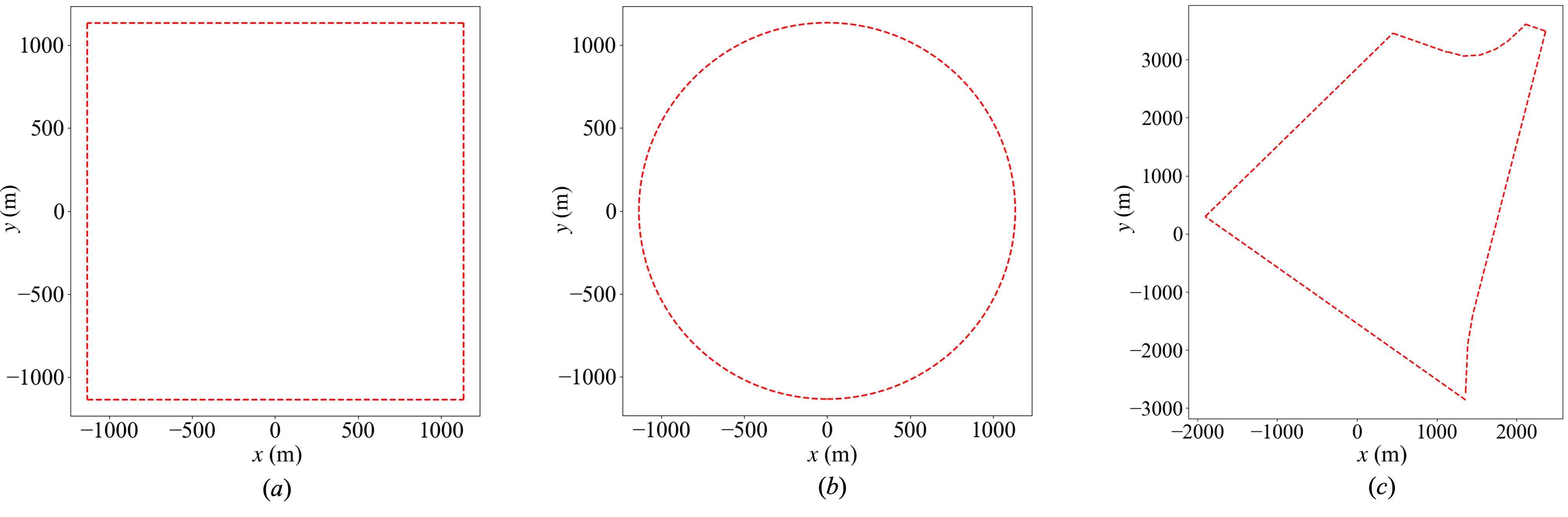}
    \caption{$(a)$ Boundary for Case I, $(b)$ boundary for Case II, $(c)$ boundary for Case III.}
    \label{fig:boundary}
\end{figure}

\begin{figure}
    \centering
    \includegraphics[width=0.9\textwidth]{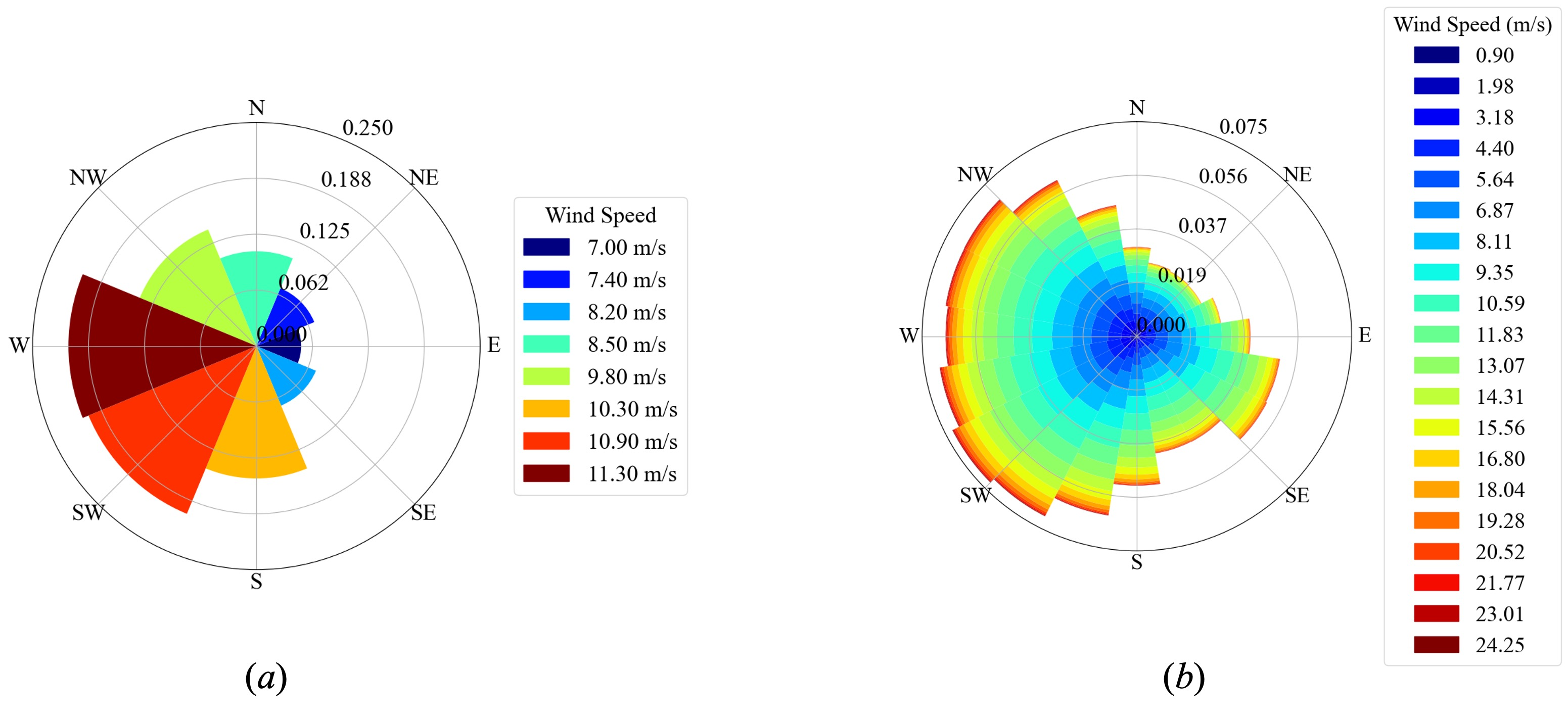}
    \caption{$(a)$ Wind rose for Case I and Case II; $(b)$ wind rose for Case III.}
    \label{fig:windrose}
\end{figure}
The NREL 5 MW reference wind turbine \citep{jonkman2009definition} is used for the three wind farm layout optimization cases.
For simplification, the wind farm area is split into $1D\times 1D$ square cells as candidate positions for wind turbines, formulating a discrete optimization problem.
The boundaries of the three wind farm cases are illustrated in detail in Figure \ref{fig:boundary}.
In Case I, a wind farm of $16$ turbines within an $18D \times 18D$ rectangular area is considered. 
Case II considers a circular area with a diameter of $18D$, in which $16$ turbines are deployed.
Case I and Case II share the same wind rose as shown in Figure \ref{fig:windrose}$(a)$.
Case III further considers a 25-turbine wind farm within a more complex boundary as shown in Figure \ref{fig:boundary}($c$), and its wind rose is shown in Figure \ref{fig:windrose}$(b)$. 
Cases II and III necessitate the use of the constrained LHS method as described in Section \ref{sec:lhs}.

\begin{figure}
    \centering
    \includegraphics[width=0.5 \textwidth]{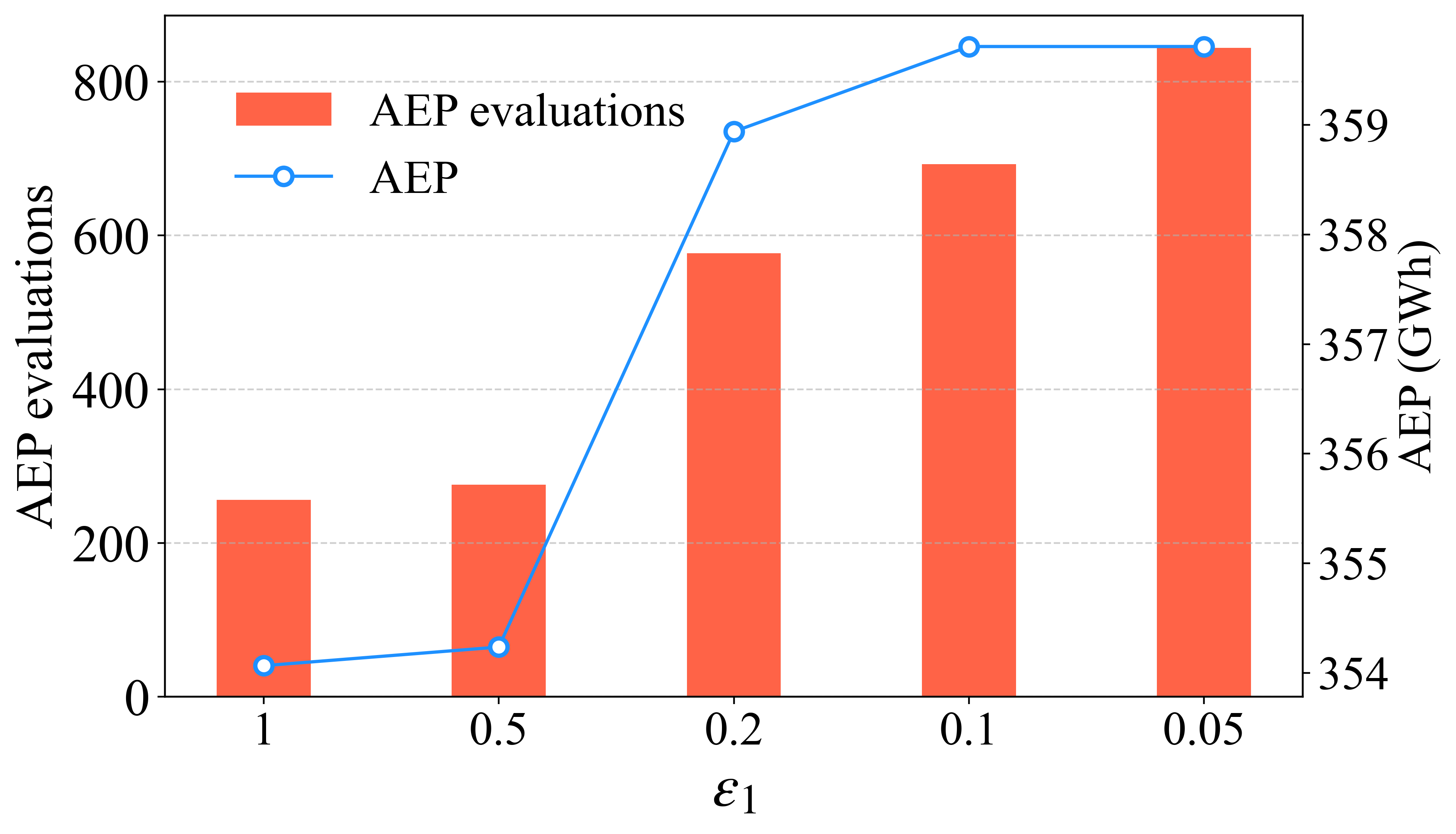}
    \caption{Sensitivity analysis for $\epsilon_1$. The AEP of optimized layouts for Case I and the corresponding number of AEP evaluations are both considered for this analysis.}
    \label{fig:epsilon1}
\end{figure}

The detailed settings for the Bayesian optimization of the three cases are illustrated in Table \ref{tab:settings}.
As the population variances limit $\epsilon_1$ significantly affects the optimized results, a sensitivity analysis based on Case I is carried out to determine its proper value in WFLO problems.
As demonstrated in Figure \ref{fig:epsilon1}, setting the $\epsilon_1$ to 0.1 yields significantly higher AEP compared to configurations where $\epsilon_1$ is set to 1, 0.5 or 0.2.
Besides, further reducing $\epsilon_1$ to 0.05 leads to an additional increase in computational cost without yielding any improvement in the results.
Therefore, $\epsilon_1$ is set as $0.1$ in this study.
The details of case setup are provided in the github repositories\footnote{\url{https://github.com/sjtuwzf/WFLBO}}.

\begin{table}
  \centering
  \caption{Detailed settings for the Bayesian optimization of three low-fidelity WFLO cases}
  \setstretch{1.2}
  \begin{tabular}{lccc}
    \toprule
    \textbf{Category} & \textbf{Case I} & \textbf{Case II} & \textbf{Case III} \\
    \midrule
    \textbf{Turbines number} & 16 & 16 & 25 \\
    \textbf{Boundary} & Figure \ref{fig:boundary}($a$) & Figure \ref{fig:boundary}($b$) & Figure \ref{fig:boundary}($c$) \\
    \textbf{Wind rose} & Figure \ref{fig:windrose}($a$) & Figure \ref{fig:windrose}($a$) & Figure \ref{fig:windrose}($b$) \\
    \midrule
    \textbf{Initial sampling} & Standard LHS, size 160 & Constrained sampling, size 160 & Constrained sampling, size 250 \\
    \midrule
    \multirow{6}{*}{\textbf{Optimization}}
    & \multicolumn{3}{c}{Optimizer: DE} \\
    & \multicolumn{3}{c}{Population size: 70} \\
    & \multicolumn{3}{c}{Mutation: (0.5,1.2)} \\
    & \multicolumn{3}{c}{Recombination: 0.7} \\
    & \multicolumn{3}{c}{Strategy: 'best1bin' for Kriging emulator, 'best2bin' for MES} \\
    & \multicolumn{3}{c}{Tolerance: $1\times 10^{-3}$ for Kriging emulator, $1\times 10^{-15}$ for MES} \\
    \midrule
    \multirow{3}{*}{\textbf{Termination criteria}}
    & \multicolumn{3}{c}{Population size: 20} \\
    & \multicolumn{3}{c}{Maximum AEP evaluations N: 1000} \\
    & \multicolumn{3}{c}{Population variances limit $\epsilon_1$: 0.1} \\
    & \multicolumn{3}{c}{MES limit $\epsilon_2$: 0.0001} \\
    \bottomrule
  \end{tabular}%
  \label{tab:settings}%
\end{table}

\subsection{Comparison with direct optimization approach}
\label{sec:direct}
In this section, optimized layouts of the three wind farm cases are obtained from the Bayesian optimization technique and three popular heuristic algorithms, i.e., differential evolution algorithm (DE), genetic algorithm (GA), and simulated annealing algorithm (SA). 
Similar to Section \ref{sec:experiments}, the results of WFLO are also influenced by the initialization (initial training samples for BO and initial population for the direct heuristic optimization), 10 different initializations are implemented for the current study.

\begin{table}
    \centering
    \captionsetup[table]{labelfont=bf,font={sf,small},singlelinecheck=off,labelsep=space}
    \captionof{table}{AEP of optimized layouts using different optimization approach and the corresponding number of AEP evaluations for 3 cases}
    \setstretch{1.5}
    \begin{tabularx}{\columnwidth}{@{\extracolsep{\fill}}*{7}{c}}
         \hline
           \multirow{2}{*}& \multicolumn{2}{c}{Case I}& \multicolumn{2}{c}{Case II}& \multicolumn{2}{c}{Case III}  \\
         \cline{2-7}
           &AEP evaluations &AEP (GWh) &AEP evaluations &AEP (GWh) &AEP evaluations &AEP (GWh) \\
         \hline
           Switch-AF       &968 &359.311 & 1160 &352.094 &1169 &559.456 \\
           DE    &10000 &339.606 &10000 &324.137 &10000 &538.934 \\
           GA  &10000 &351.744 &10000 &334.172 &10000 &544.906 \\
           SA  &10000 &364.911 &10000 &355.595 &10000 &564.812 \\
         \hline
    \end{tabularx}
    \label{tab:results}
\end{table}

Table \ref{tab:results} presents the average results of the 10 optimization runs for Case I, Case II, and Case III.
Notably, the proposed Bayesian Optimization (BO) approach required only 995, 1095, and 1169 Floris-based AEP evaluations to obtain optimized layouts for Case I, Case II, and Case III, respectively. 
It is also noted that the three heuristic algorithms fail to converge before they are terminated even if their optimized AEPs stopped improving after about $3000 - 5000$ Floris-based AEP evaluations.
This is attributed to the fact that the efficiency and effectiveness of these algorithms are highly dependent on their parameters, which is not the primary focus of this study \citep{akbaripour2013efficient}.
Nevertheless, the substantial differences between our proposed BO approach and the three heuristic algorithms demonstrate a remarkable reduction in computational effort.
In terms of performance, the optimized layouts generated by the BO approach achieved average AEP values of 98.1\%, 99.0\%, and 99.1\% relative to those obtained by SA for Case I, II, and III, respectively. 
Furthermore, the AEP values from the BO approach are substantially higher than those achieved by DE and GA.
\begin{figure}
    \centering
    \includegraphics[width=1.0 \textwidth]{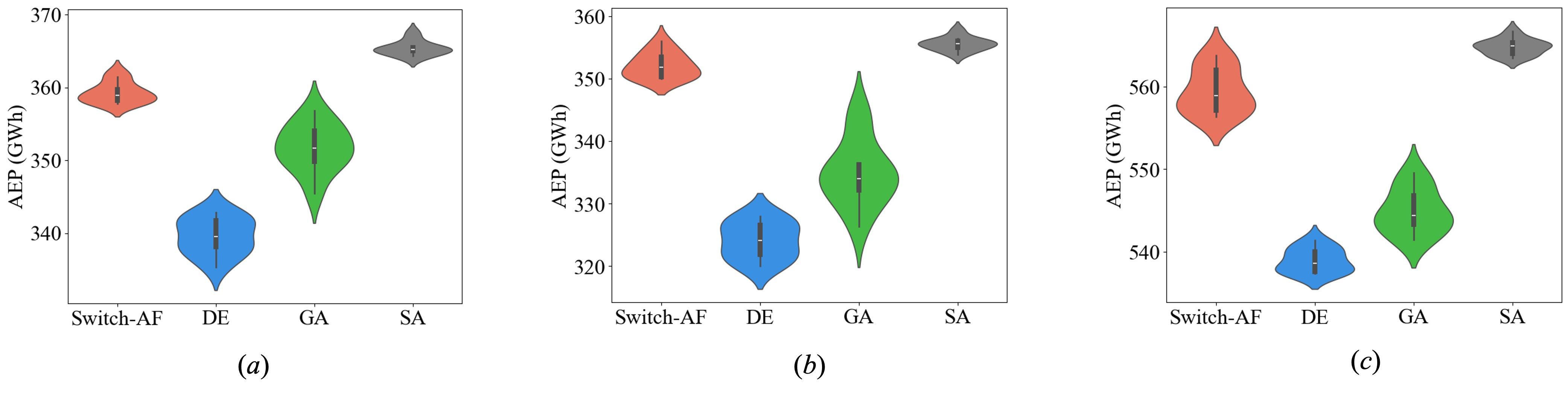}
    \caption{Violin plots for the AEPs obtained from 10 runs. $(a)$ Case I, $(b)$ Case II, and $(c)$ Case III.}
    \label{fig:violins}
\end{figure}

The violin plots in Figure \ref{fig:violins} depict the distribution of AEP values for layouts obtained from 10 random initializations.
It is observed that the AEP distributions for the BO method are more compact than those of Differential Evolution (DE) and Genetic Algorithms (GA) in Case II and Case III, though they become slightly more dispersed in Case III. 
This dispersion is primarily attributed to the complex boundary conditions in Case III, which may require significantly more computational resources to address effectively.
Despite this, the optimized layouts from all 10 runs of the BO method consistently yield higher AEP values than those of DE and GA across all three cases.

\begin{figure}
    \centering
    \includegraphics[width=1.0 \textwidth]{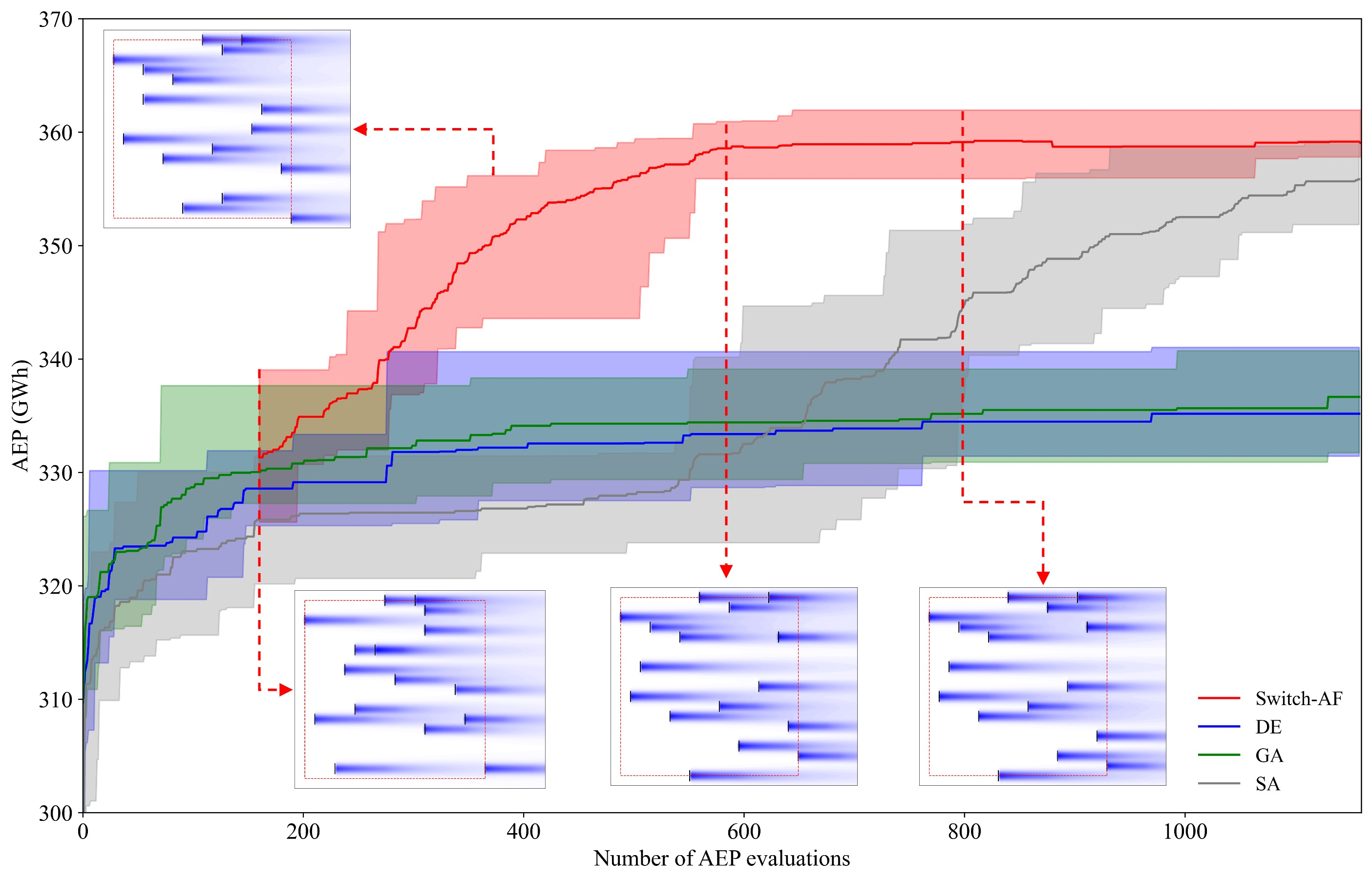}
    \caption{Evolutions of optimized AEPs obtained from switch-AF-based BO and three heuristic algorithms for Case I. The solid lines represent the average optimized AEP with 10 random initialization during the optimization; the boundaries of the shaded area represent the maximum and minimum optimized AEP with 10 random initialization during the optimization. The velocity fields of selected intermediate optimized layouts along dominant wind direction are also shown.}
    \label{fig:Case1Iteration}
\end{figure}

\begin{figure}
    \centering
    \includegraphics[width=1.0 \textwidth]{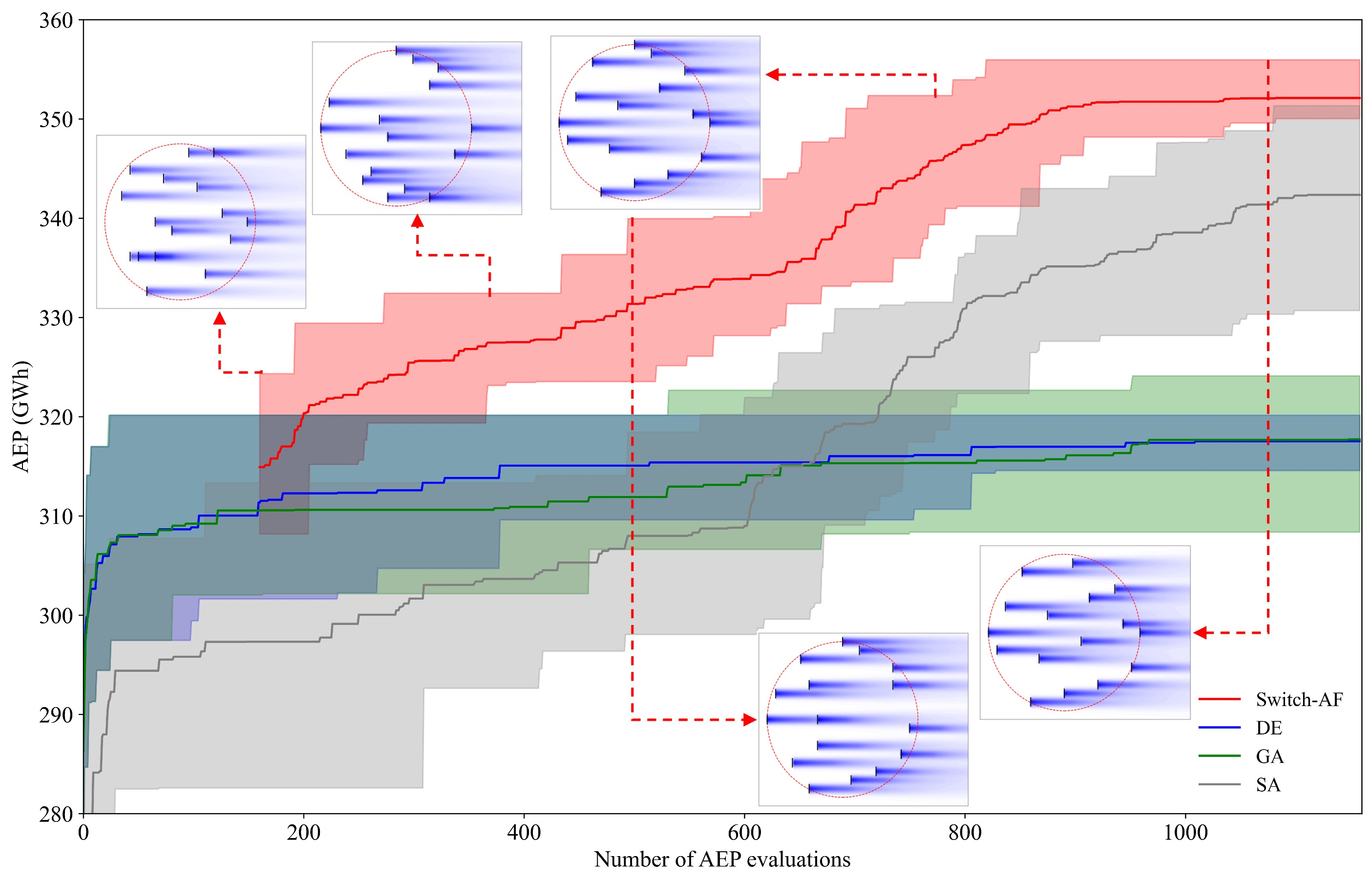}
    \caption{Evolutions of optimized AEPs obtained from switch-AF-based BO and three heuristic algorithms for Case II. }
    \label{fig:Case2Iteration}
\end{figure}

\begin{figure}
    \centering
    \includegraphics[width=1.0 \textwidth]{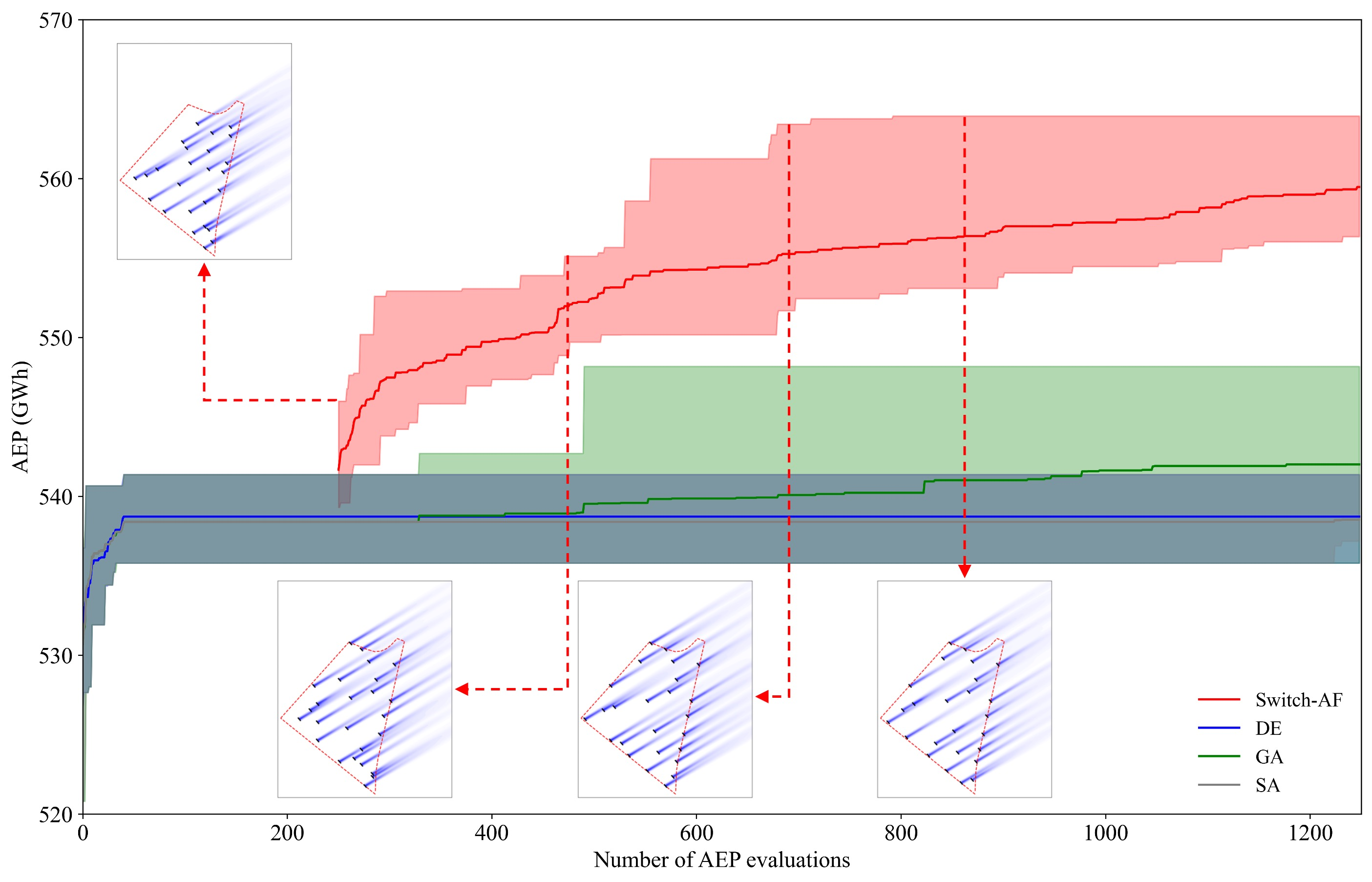}
    \caption{Evolutions of optimized AEPs obtained from switch-AF-based BO and three heuristic algorithms for Case III.}
    \label{fig:Case3Iteration}
\end{figure}

Figures \ref{fig:Case1Iteration}-\ref{fig:Case3Iteration} illustrate the progression of optimized AEP across 10 independent runs for each optimization approach in the initial stages of the optimization processes.
The results demonstrate that the optimized AEP of layouts obtained from the proposed BO framework exhibits a markedly faster rate of improvement during the initial optimization phase compared to those obtained via the three heuristic algorithms.
This is particularly true for Case III, where the optimized AEPs achieved by the three heuristic algorithms exhibit minimal changes, primarily due to the significant computational resources required by these algorithms to effectively manage the complex constraints. 
From the wake velocity fields, it is observed that as the optimization evolves, the positions of the turbines are adjusted to avoid the unfavorable effects from the wakes of the upstream turbines. 
The above discussions demonstrate that the proposed BO approach efficiently finds high-quality solutions under strict computational limits, making it ideal for resource-constrained scenarios such as CFD-based WFLO problems.

\section{CFD-based Bayesian optimization of wind farm layout}
\label{sec:high}
Building upon the previous analysis, this section incorporates the CFD simulations into the BO framework to carry out high-fidelity WFLO.
Despite the considerable cost reductions offered by the proposed approach, the inherent computational expense associated with conducting CFD simulations remains substantial. 
Thus, a single wind farm case is examined in this paper to demonstrate the applicability of the proposed framework on high-fidelity WFLO problems.

\begin{figure}
    \centering
    \includegraphics[width=0.8\linewidth]{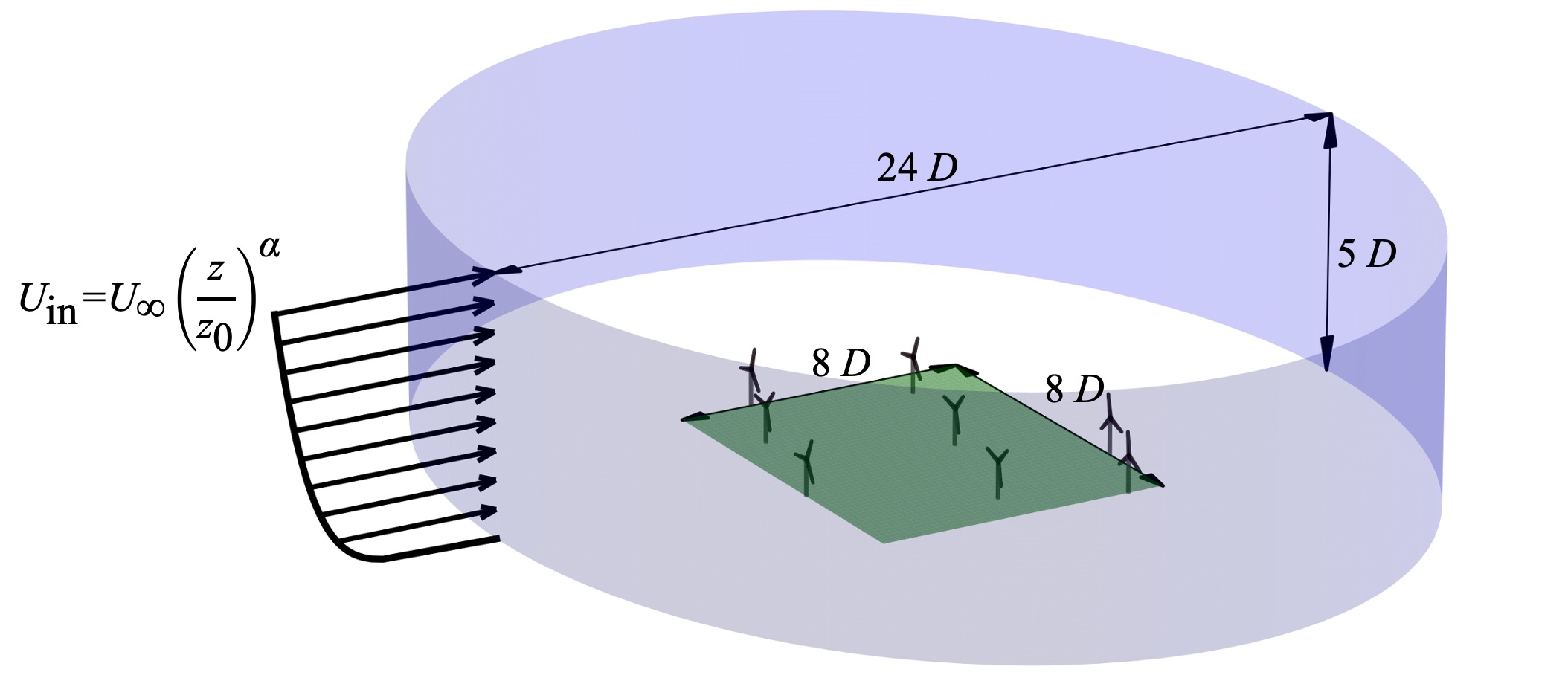}
    \caption{Setups for wind farm CFD simulations.}
    \label{fig:cfdDomain}
\end{figure}

In this case, an 8-turbine wind farm is optimized within an $8D\times 8D$ rectangular area with the same wind rose shown in Figure \ref{fig:windrose}$(a)$.
A circular computational domain encompassing the wind farm area is adopted to account for the inlet wind velocity from different directions specified by the wind rose, as shown in Figure \ref{fig:cfdDomain}.
The diameter of the circular computational domain is set as 3 times the length of the wind farm rectangle side.
In the vertical direction, the height of the computational domain is set as 5$D$.
The inlet velocity is specified as $U_{in} = U_{\infty}(z/z_0)^{\alpha}$, where $z_0 = 90$ m is the height of the turbine hub, $U_\infty$ is the reference velocity, and $\alpha = 0.14$ is the shear rate of the velocity profile.
A reference pressure $p_\infty = 0$ with zero gradient condition is set for the outlet velocity.
The mesh of the computational domain is created using OpenFOAM's native blockMesh tool.
The detailed case setup and mesh dependency test can be found in the authors' previous work \citep{wang2024optimization}.
Each CFD simulation is conducted with 60 processors using a workstation with AMD EPYC 7302 CPU, taking around 4 minutes to obtain a converged solution.

\begin{table}
    \caption{Comparison of the performance of the optimization (number of AEP evaluations and optimized AEP) between the BO using switch-AF and MSP AF, and direct optimization using GA. The results for MSP and GA are obtained from the first author's previous work \citep{wang2024optimization}, in which the same WFLO problem was considered.}
    \setstretch{1.5}
    \centering
    \begin{tabular}{ccc}
         \toprule
            & number of AEP evaluations & optimized AEP (GWh) \\
         \midrule
           switch-AF  & 301 & 127.514 \\
           MSP        & 400 & 120.374 \\
           GA         & 793 & 126.707 \\
         \bottomrule
    \end{tabular}
    \label{tab:high_results}
\end{table}

For this study, 80 different layouts are obtained using the standard LHS as the initial samples for the training of the Kriging model.
As shown in Table \ref{tab:high_results}, the switch-AF-based BO achieves an optimal solution (AEP=127.514 MW) using 221 AEP evaluations.
Since each AEP evaluation includes 8 wind directions as specified in the wind rose, the total number of CFD simulations (including the 80 training samples) amounts to 2408.
This result demonstrates significant improvement (in terms of both efficiency and optimized APE) from the authors' previous study \citep{wang2024optimization}, in which the same WFLO configuration was studied using MSP-based BO and GA with another software \citep{adams2020dakota}.

\begin{figure}
    \centering
    \includegraphics[width=1.0 \textwidth]{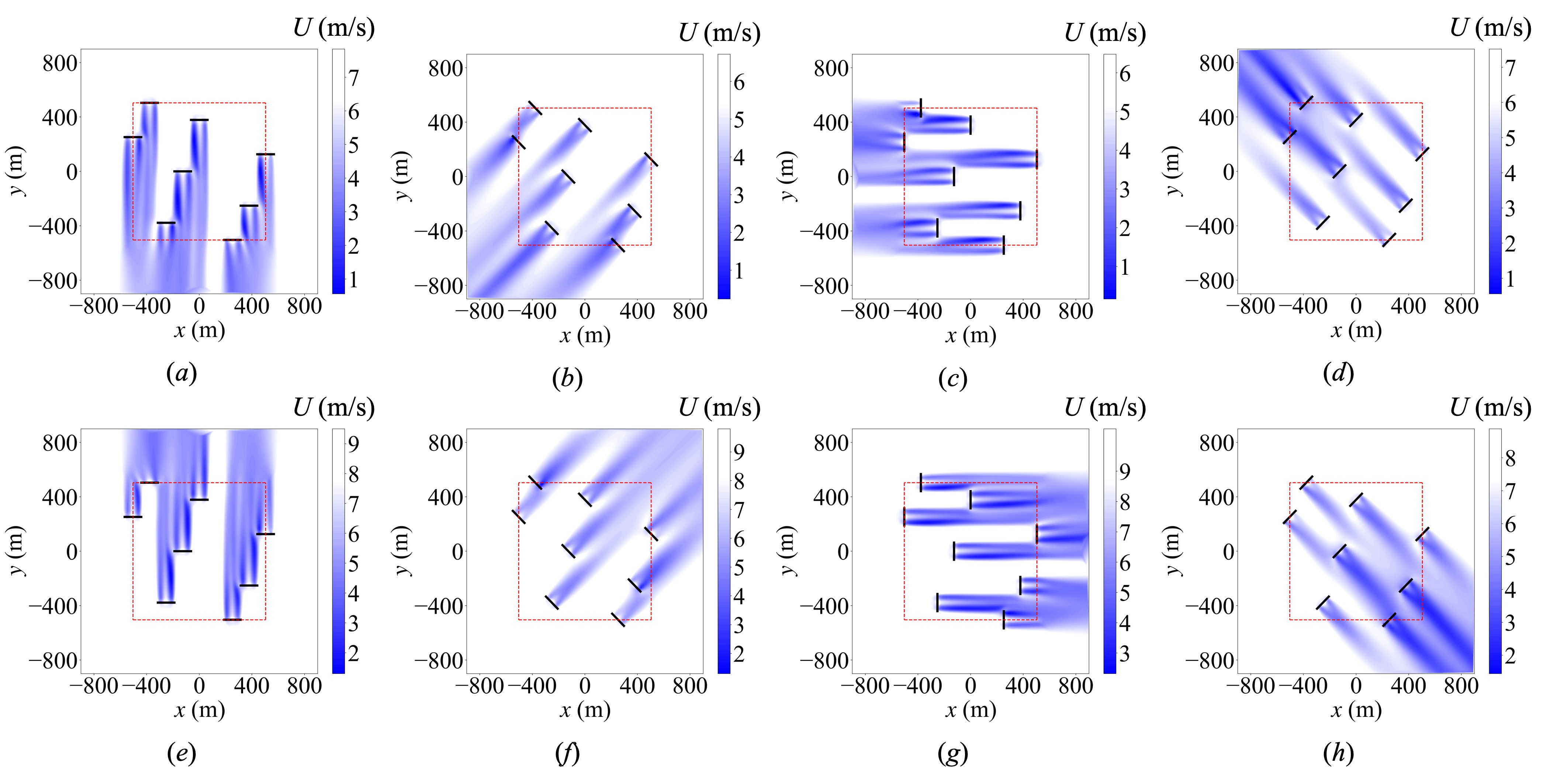}
    \caption{Velocity magnitude at rotor hub height of optimized layouts in the high-fidelity WFLO case. ($a$) northern, ($b$) northeastern ($c$) eastern, ($d$) southeastern, ($e$) southern, ($f$) southwestern, ($g$) western (dominant wind direction) and ($f$) northwestern inflow. The red dashed line represents the wind farm boundary.}
    \label{fig:CFDWwake}
\end{figure}

The velocity fields at the hub height of optimized layouts in Figure \ref{fig:CFDWwake}.
Compared with the flow fields generated by low-fidelity analytical wake models as shown in Figures \ref{fig:Case1Iteration}-\ref{fig:Case3Iteration}, the ADM-based CFD simulations capture more intrinsic features such as higher velocity near the blade root region and wake asymmetry due to interactions among the turbine cluster.
Along most of the wind directions, the wind turbines are strategically positioned so that the wakes of the upstream ones barely affect their downstream peers.
The successful application of the switch-AF-based Bayesian optimization with CFD simulations demonstrates its potential to enable high-fidelity wind farm layout design.

\section{Conclusions}
\label{sec:conclusion}

This paper developed a novel adaptive switch strategy for the acquisition functions in the Bayesian optimization (BO) framework, aiming at facilitating high-fidelity wind farm layout optimization.
The Kriging model is used as the surrogate for efficient AEP prediction during optimization, and an adaptive switch acquisition function strategy between Max-value Entropy Search (MES) and Max Surrogate Prediction (MSP) is proposed to achieve a balance between exploration and exploitation in the infilling criterion in Bayesian optimization.
The Kriging's AEP emulator, which is also effectively the MSP acquisition function, is first integrated with the differential evolution (DE) algorithm, and the intermediate design during the optimization process is put back into the training samples to update the Kriging model.
Once optimization using the MSP AF converges, the MES AF assumes control to refine the search for global optima. 

The proposed switch-AF-based BO framework is first validated on 4D and 10D Ackley functions, demonstrating its enhanced optimization efficiency over the conventional use of either MSP or MES as standalone acquisition functions. 
Subsequently, the framework is evaluated using computationally efficient Gaussian wake models for WFLO across three wind farm cases with varying boundaries and wind distributions. 
The results reveal that the switch-AF-based BO framework outperforms traditional heuristic algorithms (DE, GA, and SA) in both computational efficiency and optimization effectiveness, achieving near-optimal AEP with significantly fewer evaluations. 
Moreover, the integration of CFD simulations into the BO framework highlights its potential for high-fidelity WFLO, providing detailed insights into wake interactions and turbine placement.

While the current study demonstrates that Bayesian optimization can reduce the number of CFD simulations from thousands to hundreds, challenges remain for larger wind farms. The computational cost for hundreds of turbines remains prohibitively high, even with reduced simulations. Future work should focus on multi-fidelity Kriging models, combining low-fidelity data with limited high-fidelity data, to extend the framework's applicability. Additionally, the curse of dimensionality in high-dimensional problems necessitates dimensionality reduction techniques to improve the efficiency of Gaussian process regression for large-scale WFLO.

In conclusion, this paper introduces an adaptive BO framework that significantly reduces computational costs while maintaining high-fidelity accuracy for wind farm layout optimizations. 
By balancing exploration and exploitation through a novel switch strategy, the framework offers a robust and efficient solution for optimizing wind farm layouts, paving the way for more sustainable and effective wind energy projects.

\section*{Acknowledgments}{The financial supports from the National Key R\&D Program of China (No. 2023YFE0120000), Guangdong Basic and Applied Basic Research Foundation (No. 2023A1515240054), National Natural Science Foundation of China (No. 42076210), and Program for Intergovernmental International S\&T Cooperation Projects of Shanghai Municipality, China (No. 22160710200) are gratefully acknowledged. }\label{sec:acknowledgments}

\section*{Declaration of interests}{

The authors report no conflict of interest.

}

\bibliographystyle{unsrtnat}
\bibliography{reference}

\begin{thebibliography}{46}
\providecommand{\natexlab}[1]{#1}
\providecommand{\url}[1]{\texttt{#1}}
\expandafter\ifx\csname urlstyle\endcsname\relax
  \providecommand{\doi}[1]{doi: #1}\else
  \providecommand{\doi}{doi: \begingroup \urlstyle{rm}\Url}\fi

\bibitem[Spiru and Simona(2024)]{spiru2024wind}
Paraschiv Spiru and Paraschiv~Lizica Simona.
\newblock Wind energy resource assessment and wind turbine selection analysis
  for sustainable energy production.
\newblock \emph{Scientific Reports}, 14\penalty0 (1):\penalty0 10708, 2024.

\bibitem[Barthelmie et~al.(2009)Barthelmie, Hansen, Frandsen, Rathmann,
  Schepers, Schlez, Phillips, Rados, Zervos, Politis,
  et~al.]{barthelmie2009modelling}
Rebecca~Jane Barthelmie, K~Hansen, Sten~Tron{\ae}s Frandsen, Ole Rathmann,
  JG~Schepers, W~Schlez, J~Phillips, K~Rados, A~Zervos, ESa Politis, et~al.
\newblock Modelling and measuring flow and wind turbine wakes in large wind
  farms offshore.
\newblock \emph{Wind Energy}, 12\penalty0 (5):\penalty0 431--444, 2009.

\bibitem[Ero{\u{g}}lu and Se{\c{c}}kiner(2012)]{erouglu2012design}
Yunus Ero{\u{g}}lu and Serap~Ulusam Se{\c{c}}kiner.
\newblock Design of wind farm layout using ant colony algorithm.
\newblock \emph{Renewable Energy}, 44:\penalty0 53--62, 2012.

\bibitem[Stanley and Ning(2019)]{stanley2019massive}
Andrew~PJ Stanley and Andrew Ning.
\newblock Massive simplification of the wind farm layout optimization problem.
\newblock \emph{Wind Energy Science}, 4\penalty0 (4):\penalty0 663--676, 2019.

\bibitem[Antonini et~al.(2020)Antonini, Romero, and Amon]{antonini2020optimal}
Enrico~GA Antonini, David~A Romero, and Cristina~H Amon.
\newblock Optimal design of wind farms in complex terrains using computational
  fluid dynamics and adjoint methods.
\newblock \emph{Applied Energy}, 261:\penalty0 114426, 2020.

\bibitem[Cruz and Carmo(2020)]{cruz2020wind}
Lu{\'\i}s Eduardo~Boni Cruz and Bruno~Souza Carmo.
\newblock Wind farm layout optimization based on {CFD} simulations.
\newblock \emph{Journal of the Brazilian Society of Mechanical Sciences and
  Engineering}, 42\penalty0 (8):\penalty0 433, 2020.

\bibitem[Reddy(2020)]{reddy2020wind}
Sohail~R Reddy.
\newblock Wind farm layout optimization ({WindFLO}): {An} advanced framework
  for fast wind farm analysis and optimization.
\newblock \emph{Applied Energy}, 269:\penalty0 115090, 2020.

\bibitem[Dong et~al.(2021)Dong, Zhang, and Zhao]{dong2021intelligent}
Hongyang Dong, Jincheng Zhang, and Xiaowei Zhao.
\newblock Intelligent wind farm control via deep reinforcement learning and
  high-fidelity simulations.
\newblock \emph{Applied Energy}, 292:\penalty0 116928, 2021.

\bibitem[Thomas et~al.(2023)Thomas, Baker, Malisani, Quaeghebeur, Sanchez
  Perez-Moreno, Jasa, Bay, Tilli, Bieniek, Robinson,
  et~al.]{thomas2023comparison}
Jared~J Thomas, Nicholas~F Baker, Paul Malisani, Erik Quaeghebeur, Sebastian
  Sanchez Perez-Moreno, John Jasa, Christopher Bay, Federico Tilli, David
  Bieniek, Nick Robinson, et~al.
\newblock A comparison of eight optimization methods applied to a wind farm
  layout optimization problem.
\newblock \emph{Wind Energy Science}, 8\penalty0 (5):\penalty0 865--891, 2023.

\bibitem[Jensen(1983)]{jensen1983note}
Niels~Otto Jensen.
\newblock \emph{A note on wind generator interaction}.
\newblock Report Risø-M-2411, Risø National Laboratory, Roskilde, Denmark,
  1983.

\bibitem[Larsen(1988)]{larsen1988simple}
Gunner~Chr Larsen.
\newblock \emph{{A Simple Wake Calculation Procedure}}.
\newblock Report Riso-M-2760, Ris{\o} National Laboratory, Roskilde, Denmark,
  1988.

\bibitem[Mart{\'\i}nez-Tossas et~al.(2015)Mart{\'\i}nez-Tossas, Churchfield,
  and Leonardi]{martinez2015large}
Luis~A Mart{\'\i}nez-Tossas, Matthew~J Churchfield, and Stefano Leonardi.
\newblock Large eddy simulations of the flow past wind turbines: actuator line
  and disk modeling.
\newblock \emph{Wind Energy}, 18\penalty0 (6):\penalty0 1047--1060, 2015.

\bibitem[Stevens et~al.(2018)Stevens, Mart{\'\i}nez-Tossas, and
  Meneveau]{stevens2018comparison}
Richard J. A.~M. Stevens, Luis~A Mart{\'\i}nez-Tossas, and Charles Meneveau.
\newblock Comparison of wind farm large eddy simulations using actuator disk
  and actuator line models with wind tunnel experiments.
\newblock \emph{Renewable Energy}, 116:\penalty0 470--478, 2018.

\bibitem[Tu et~al.(2023)Tu, Zhang, Han, Zhou, and Bilgen]{tu2023aerodynamic}
Yu~Tu, Kai Zhang, Zhaolong Han, Dai Zhou, and Onur Bilgen.
\newblock Aerodynamic characterization of two tandem wind turbines under yaw
  misalignment control using actuator line model.
\newblock \emph{Ocean Engineering}, 281:\penalty0 114992, 2023.

\bibitem[Wang et~al.(2024)Wang, Tu, Zhang, Han, Cao, and
  Zhou]{wang2024optimization}
Zhenfan Wang, Yu~Tu, Kai Zhang, Zhaolong Han, Yong Cao, and Dai Zhou.
\newblock An optimization framework for wind farm layout design using
  {CFD}-based {Kriging} model.
\newblock \emph{Ocean Engineering}, 293:\penalty0 116644, 2024.

\bibitem[Seeger(2004)]{gaussian2004}
MATTHIAS Seeger.
\newblock Gaussian processes for machine learning.
\newblock \emph{International Journal of Neural Systems}, 14\penalty0
  (02):\penalty0 69--106, 2004.

\bibitem[Frazier(2018)]{frazier2018tutorial}
Peter~I Frazier.
\newblock A tutorial on {Bayesian} optimization.
\newblock \emph{arXiv preprint arXiv:1807.02811}, 2018.

\bibitem[Kushner(1964)]{kushner1964new}
H.~J. Kushner.
\newblock A new method of locating the maximum point of an arbitrary multipeak
  curve in the presence of noise.
\newblock \emph{Journal of Basic Engineering}, 86\penalty0 (1):\penalty0
  97--106, 03 1964.
\newblock ISSN 0021-9223.

\bibitem[Mockus(1974)]{mockus1974bayesian}
Jonas Mockus.
\newblock On {Bayesian} methods for seeking the extremum.
\newblock In \emph{Proceedings of the IFIP Technical Conference}, pages
  400--404, 1974.

\bibitem[Shahriari et~al.(2016)Shahriari, Swersky, Wang, Adams, and
  de~Freitas]{shahriari2015taking}
Bobak Shahriari, Kevin Swersky, Ziyu Wang, Ryan~P. Adams, and Nando de~Freitas.
\newblock Taking the human out of the loop: A review of {Bayesian}
  optimization.
\newblock \emph{Proceedings of the IEEE}, 104\penalty0 (1):\penalty0 148--175,
  2016.

\bibitem[Qin et~al.(2017)Qin, Klabjan, and Russo]{qin2017improving}
Chao Qin, Diego Klabjan, and Daniel Russo.
\newblock Improving the expected improvement algorithm.
\newblock \emph{Advances in Neural Information Processing Systems}, 30, 2017.

\bibitem[Hennig and Schuler(2012)]{hennig2012entropy}
Philipp Hennig and Christian~J Schuler.
\newblock Entropy search for information-efficient global optimization.
\newblock \emph{Journal of Machine Learning Research}, 13\penalty0 (6), 2012.

\bibitem[Gallager(1968)]{gallager1968information}
Robert~G Gallager.
\newblock \emph{Information Theory and Reliable Communication}, volume 588.
\newblock Wiley, 1968.

\bibitem[Hern{\'a}ndez-Lobato et~al.(2016)Hern{\'a}ndez-Lobato,
  Hernandez-Lobato, Shah, and Adams]{hernandez2016predictive}
Daniel Hern{\'a}ndez-Lobato, Jose Hernandez-Lobato, Amar Shah, and Ryan Adams.
\newblock Predictive entropy search for multi-objective bayesian optimization.
\newblock In \emph{Proceedings of The 33rd International Conference on Machine
  Learning}, pages 1492--1501, 2016.

\bibitem[Wang and Jegelka(2017)]{wang2017max}
Zi~Wang and Stefanie Jegelka.
\newblock Max-value entropy search for efficient {Bayesian} optimization.
\newblock In \emph{International Conference on Machine Learning}, pages
  3627--3635. PMLR, 2017.

\bibitem[Rehbach et~al.(2020)Rehbach, Zaefferer, Naujoks, and
  Bartz-Beielstein]{rehbach2020expected}
Frederik Rehbach, Martin Zaefferer, Boris Naujoks, and Thomas Bartz-Beielstein.
\newblock Expected improvement versus predicted value in surrogate-based
  optimization.
\newblock In \emph{Proceedings of the 2020 genetic and evolutionary computation
  conference}, pages 868--876, 2020.

\bibitem[Bempedelis et~al.(2024)Bempedelis, Gori, Wynn, Laizet, and
  Magri]{bempedelis2024data}
N.~Bempedelis, F.~Gori, A.~Wynn, S.~Laizet, and L.~Magri.
\newblock Data-driven optimisation of wind farm layout and wake steering with
  large-eddy simulations.
\newblock \emph{Wind Energy Science}, 9\penalty0 (4):\penalty0 869--882, 2024.

\bibitem[Krige(1951)]{krige1951statistical}
Daniel~G Krige.
\newblock \emph{A statistical approach to some mine valuation and allied
  problems on the {Witwatersrand}}.
\newblock PhD thesis, University of the Witwatersrand, 1951.

\bibitem[Rasmussen and Williams(2006)]{Rasmussen2006Gaussian}
Carl~Edward Rasmussen and Christopher K.~I. Williams.
\newblock \emph{{Gaussian Processes for Machine Learning}}.
\newblock The MIT Press, 2006.

\bibitem[Adams et~al.(2020)Adams, Bohnhoff, Dalbey, Ebeida, Eddy, Eldred,
  Hooper, Hough, Hu, Jakeman, et~al.]{adams2020dakota}
Brian~M Adams, William~J Bohnhoff, Keith~R Dalbey, Mohamed~S Ebeida, John~P
  Eddy, Michael~S Eldred, Russell~W Hooper, Patricia~D Hough, Kenneth~T Hu,
  John~D Jakeman, et~al.
\newblock Dakota, a multilevel parallel object-oriented framework for design
  optimization, parameter estimation, uncertainty quantification, and
  sensitivity analysis: version 6.13 user's manual.
\newblock Technical report, Sandia National Lab, Albuquerque, NM, 2020.

\bibitem[Sacks et~al.(1989)Sacks, Schiller, and Welch]{KrigingTrain}
Jerome Sacks, Susannah~B. Schiller, and William~J. Welch.
\newblock Designs for {Computer Experiments}.
\newblock \emph{Technometrics}, 31\penalty0 (1):\penalty0 41--47, 1989.
\newblock ISSN 00401706.

\bibitem[Saves et~al.(2024)Saves, Lafage, Bartoli, Diouane, Bussemaker,
  Lefebvre, Hwang, Morlier, and Martins]{saves2024smt}
Paul Saves, R{\'e}mi Lafage, Nathalie Bartoli, Youssef Diouane, Jasper
  Bussemaker, Thierry Lefebvre, John~T Hwang, Joseph Morlier, and Joaquim~RRA
  Martins.
\newblock {SMT 2.0: A Surrogate Modeling Toolbox with a focus on hierarchical
  and mixed variables Gaussian processes}.
\newblock \emph{Advances in Engineering Software}, 188:\penalty0 103571, 2024.

\bibitem[McKay et~al.(2000)McKay, Beckman, and Conover]{mckay2000comparison}
Michael~D McKay, Richard~J Beckman, and William~J Conover.
\newblock A comparison of three methods for selecting values of input variables
  in the analysis of output from a computer code.
\newblock \emph{Technometrics}, 42\penalty0 (1):\penalty0 55--61, 2000.

\bibitem[Conde~Arenzana et~al.(2021)Conde~Arenzana, L{\'o}pez-Lopera, Mouton,
  Bartoli, and Lefebvre]{arenzana2021multi}
Rub{\'e}n Conde~Arenzana, Andr{\'e}s~F. L{\'o}pez-Lopera, Sylvain Mouton,
  Nathalie Bartoli, and Thierry Lefebvre.
\newblock Multi-fidelity {Gaussian Process} model for {CFD} and wind tunnel
  data fusion.
\newblock In \emph{ECCOMAS AeroBest}, 2021.

\bibitem[Gray(2011)]{gray2011entropy}
Robert~M Gray.
\newblock \emph{Entropy and Information Theory}.
\newblock Springer Science \& Business Media, 2011.

\bibitem[Virtanen et~al.(2020)Virtanen, Gommers, Oliphant, Haberland, Reddy,
  Cournapeau, Burovski, Peterson, Weckesser, Bright, et~al.]{2020SciPy-NMeth}
Pauli Virtanen, Ralf Gommers, Travis~E Oliphant, Matt Haberland, Tyler Reddy,
  David Cournapeau, Evgeni Burovski, Pearu Peterson, Warren Weckesser, Jonathan
  Bright, et~al.
\newblock {SciPy 1.0: fundamental algorithms for scientific computing in
  Python}.
\newblock \emph{Nature Methods}, 17\penalty0 (3):\penalty0 261--272, 2020.

\bibitem[Motiian and Soltanian-Zadeh(2011)]{motiian2011}
Saeed Motiian and Hamid Soltanian-Zadeh.
\newblock Improved particle swarm optimization and applications to hidden
  markov model and ackley function.
\newblock In \emph{IEEE International Conference on Computational Intelligence
  for Measurement Systems and Applications Proceedings}, pages 1--4, 2011.

\bibitem[King et~al.(2021)King, Fleming, King, Mart{\'\i}nez-Tossas, Bay,
  Mudafort, and Simley]{king2021control}
Jennifer King, Paul Fleming, Ryan King, Luis~A Mart{\'\i}nez-Tossas,
  Christopher~J Bay, Rafael Mudafort, and Eric Simley.
\newblock Control-oriented model for secondary effects of wake steering.
\newblock \emph{Wind Energy Science}, 6\penalty0 (3):\penalty0 701--714, 2021.

\bibitem[NREL(2023)]{Floris}
NREL.
\newblock {FLORIS}, 12 2023.
\newblock URL \url{https://github.com/nrel/floris}.

\bibitem[Launder and Spalding(1974)]{Launder}
Brian Launder and D.B. Spalding.
\newblock The numerical computation of turbulent flow computer methods.
\newblock \emph{Computer Methods in Applied Mechanics and Engineering},
  3:\penalty0 269--289, 03 1974.

\bibitem[Burton et~al.(2011)Burton, Jenkins, Sharpe, and
  Bossanyi]{burton2011wind}
Tony Burton, Nick Jenkins, David Sharpe, and Ervin Bossanyi.
\newblock \emph{Wind Energy Handbook}.
\newblock Wiley, 2011.

\bibitem[Goldstein(1929)]{goldstein1929vortex}
Sydney Goldstein.
\newblock On the vortex theory of screw propellers.
\newblock \emph{Proceedings of the Royal Society of London. Series A,
  Containing Papers of a Mathematical and Physical Character}, 123\penalty0
  (792):\penalty0 440--465, 1929.

\bibitem[Richmond et~al.(2019)Richmond, Antoniadis, Wang, Kolios, Al-Sanad, and
  Parol]{richmond2019evaluation}
Mark Richmond, A~Antoniadis, Lin Wang, Athanasios Kolios, S~Al-Sanad, and
  Jafarali Parol.
\newblock Evaluation of an offshore wind farm computational fluid dynamics
  model against operational site data.
\newblock \emph{Ocean Engineering}, 193:\penalty0 106579, 2019.

\bibitem[Weller et~al.(1998)Weller, Tabor, Jasak, and
  Fureby]{Weller1998OpenFOAM}
H.~G. Weller, G.~Tabor, H.~Jasak, and C.~Fureby.
\newblock A tensorial approach to computational continuum mechanics using
  object-oriented techniques.
\newblock \emph{Computer in Physics}, 12\penalty0 (6):\penalty0 620--631, 11
  1998.
\newblock ISSN 0894-1866.

\bibitem[Jonkman(2009)]{jonkman2009definition}
J~Jonkman.
\newblock {Definition of a 5-MW Reference Wind Turbine for Offshore System
  Development}.
\newblock \emph{National Renewable Energy Laboratory}, 2009.

\bibitem[Akbaripour et~al.(2013)Akbaripour, Masehian, and
  Akbaripour]{akbaripour2013efficient}
Hossein Akbaripour, E.~Masehian, and Hossein Akbaripour.
\newblock Efficient and robust parameter tuning for heuristic algorithms.
\newblock \emph{International Journal of Industrial Engineering \& Production
  Research}, 24:\penalty0 143--150, 2013.

\end{thebibliography}

\end{document}